\documentclass[twocolumn]{autart}    

\usepackage{graphicx}          
\usepackage{amsmath,amssymb,mathtools,amsfonts}
\usepackage{hyperref}
\usepackage{multirow}
\usepackage{color}
\usepackage{psfrag}
\usepackage{empheq}
\usepackage{tabularx}
\usepackage{bm} 
\usepackage{url}
\usepackage{float}
\usepackage{comment}

\usepackage[short]{optidef}
\usepackage[shortlabels]{enumitem}

\usepackage{tikz}
\usepackage{pgfplots}
\usetikzlibrary{arrows.meta}
\newlength\fwidth
\newlength\fheight
\usetikzlibrary{arrows}

\definecolor{offwhite}{RGB}{249,242,215}
\definecolor{foreground}{RGB}{0,0,0}
\definecolor{background}{RGB}{255,255,255}
\definecolor{title}{RGB}{0,0,0}
\definecolor{gray}{RGB}{155,155,155}
\definecolor{subtitle}{RGB}{155,155,155}
\definecolor{hilight}{RGB}{102,255,204}
\definecolor{vhilight}{RGB}{255,111,207}
\definecolor{lolight}{RGB}{155,155,155}
\definecolor{blueg}{RGB}{114,143,178}
\definecolor{bluet}{RGB}{100,121,178}
\definecolor{redt}{RGB}{255,86,73}
\definecolor{greent}{RGB}{136,255,134}
\definecolor{dgreent}{RGB}{75,212,69}
\definecolor{lgray}{RGB}{190,190,190}
\definecolor{llgray}{RGB}{220,220,220}

\newcommand{\R}{{\mathbb{R}}}
\newcommand{\dd}{{\mathrm{d}}}
\newcommand{\dt}{{\frac{\dd}{\dd t}}}

\newcommand{\ts}{{t_{\mathrm{s}}}}
\newcommand{\tf}{{t_{\mathrm{f}}}}

\newcommand{\qs}{{q_{\mathrm{s}}}}
\newcommand{\vs}{{v_{\mathrm{s}}}}

\newcommand{\taus}{{\tau_{\mathrm{s}}}}
\newcommand{\tauf}{{\tau_{\mathrm{f}}}}
\newcommand{\taue}{{\tau_{\mathrm{e}}}}
\newcommand{\taur}{{\tau_{\mathrm{r}}}}
\newcommand{\taujump}{{\tau_{\mathrm{jump}}}}

\newcommand{\zn}{\lambda_{\mathrm{n}}}
\newcommand{\zt}{\lambda_{\mathrm{t}}}
\newcommand{\zttilde}{\tilde{\lambda}_{\mathrm{t}}}
\newcommand{\vt}{v_{\mathrm{t}}}
\newcommand{\Tan}{B(q)}

\newcommand{\nt}{n_{\mathrm{t}}}
\newcommand{\epst}{\epsilon_{\mathrm{t}}}

\newcommand{\anormal}{a_{\mathrm{n}}}
\newcommand{\atangent}{a_{\mathrm{t}}}

\newcommand{\sign}{\mathrm{sign}}
\newcommand{\Sign}{\mathrm{Sign}}
\newcommand{\Nctrl}{{N}}
\newcommand{\Tctrl}{T}
\newcommand{\lambdap}{\lambda^{\mathrm{p}}}
\newcommand{\lambdan}{\lambda^{\mathrm{n}}}

\date{}

\newcommand{\fslip}{{{f_{\mathrm{Slip}}}}}
\newcommand{\fstick}{{{f_{\mathrm{Stick}}}}}

\newtheorem{theorem}{Theorem}
\newtheorem{example}{Example}
\newtheorem{remark}[theorem]{Remark}
\newtheorem{definition}[theorem]{Definition}

\newtheorem{proposition}[theorem]{Proposition}

\journal{Automatica}

\begin{document}

\begin{frontmatter}
\runtitle{Time-Freezing for CLS}

\title{The Time-Freezing Reformulation for Numerical Optimal Control of Complementarity Lagrangian Systems with State Jumps.\thanksref{footnoteinfo}} 

\thanks[footnoteinfo]{This research was supported by the German Federal Ministry of Education and Research (BMBF) via the funded Kopernikus project: SynErgie (03SFK3U0), DFG via Research Unit FOR 2401 and project 424107692 and by the EU via ELO-X 953348. Corresponding author A.~Nurkanovi\'c.\\
\textit{Email address}: \texttt{armin.nurkanovic@imtek.uni-freiburg.de} (Armin Nurkanovi\'c),  
\texttt{sebastian.albrecht@siemens.com} (Sebastian Albrecht),  
\texttt{bernard.brogliato@inria.fr} (Bernard Brogliato),  
\texttt{moritz.diehl@imtek.uni-freiburg.de} (Moritz Diehl)  
}

\author[Freiburg]{Armin Nurkanovi\'c},
\author[Munich]{Sebastian Albrecht},
\author[Grenoble]{Bernard Brogliato},
\author[Freiburg,Freiburg1]{Moritz Diehl},

\address[Freiburg]{Department of Microsystems Engineering (IMTEK), University of Freiburg, Germany}        
\address[Freiburg1]{Department of Mathematics, University of Freiburg, Germany}        
\address[Munich]{Siemens Technology, Munich, Germany}  
\address[Grenoble]{Univ. Grenoble Alpes, INRIA, CNRS, Grenoble INP, LJK, Grenoble, France}        
\vspace{-0.05cm}

\begin{abstract}
This paper introduces a novel time-freezing reformulation and numerical methods for optimal control of complementarity Lagrangian systems (CLS) with state jumps.
We cover the difficult case when the system evolves on the boundary of the dynamic's feasible set after the state jump.
In nonsmooth mechanics, this corresponds to inelastic impacts. 
The main idea of the time-freezing reformulation is to introduce a clock state and an auxiliary dynamical system whose trajectory endpoints satisfy the state jump law.
When the auxiliary system is active, the clock state is not evolving, hence by taking only the parts of the trajectory when the clock state was active, we can recover the original solution. 
The resulting time-freezing system is a Filippov system that has jump discontinuities only in the first time derivative instead of the trajectory itself. 
This enables one to use the recently proposed Finite Elements with Switch Detection \cite{Nurkanovic2022}, which makes high accuracy numerical optimal control of CLS with impacts and friction possible. 
We detail how to recover the solution of the original system and show how to select appropriate auxiliary dynamics.
The theoretical findings are illustrated on a nontrivial numerical optimal control example of a hopping one-legged robot.
\end{abstract}
\end{frontmatter}
\section{Introduction}\label{sec:introduction}
Complementarity Lagrangian Systems (CLS) model the dynamics of rigid bodies with friction and impact.
They are indispensable in modern robotic control applications, as any complex task requires exploiting contacts and friction~\cite{Brogliato2016,Halm2021,Howell2022,Posa2014,Stewart2000,Tassa2012}.
Although many mature and efficient simulation methods \cite{Acary2008,Stewart1996} and software \cite{Howell2022} exist, solving optimal control problems remains very difficult and is often heuristic-based.
In this paper, we introduce a time-freezing reformulation that transforms the CLS with state jumps into a Piecewise Smooth System (PSS), which is treated as a Filippov system \cite{Filippov1988}.
The resulting system has discontinuities only in its vector field, but not in the solution itself anymore. 
Our main motivation is to formulate optimal control problems that we can solve with high accuracy and to avoid convergence to spurious solutions, as we discuss below in more detail.

From the theoretical side, it is desirable to know whether solutions to optimal control problems involving CLS or Filippov systems exist. 
Currently, to the best of the authors' knowledge, there are no general results for these problems, but only for some special cases. 
This is an intriguing and challenging problem but beyond the scope of this paper.  
We will assume that optimal solutions do exist.
Under this assumption, Pontryagin-type conditions for dynamic complementarity systems (with absolutely continuous solutions) are provided in \cite{Guo2016,Vieira2019}. 
Necessary optimality conditions for a controlled first-order sweeping process \cite{Brogliato2020} are derived in \cite{Colombo2020}. 
Existence results for differential inclusions with Lipschitz properties are derived in \cite{Clarke1976,Mordukhovich2016} and with upper semi-continuity properties in~\cite{Cellina2003}.

From the computational side, time-stepping methods for CLS  and differential inclusions have in general first-order accuracy \cite{Acary2008,Stewart2000}. 
Despite this, they are often used in direct optimal control approaches \cite{Nurkanovic2020,Posa2014,Tassa2012,Vieira2019}. 
However, it was first noticed by Stewart and Anitescu \cite{Stewart2010} that for ODE with a discontinuous right-hand side (r.h.s.), contrary to smooth ODE, the numerical sensitivities are wrong independent of the step size. 
This often impairs the progress of the optimizer and results in spurious local minima~\cite{Nurkanovic2020}.
Even smoothing provides serious pitfalls besides the introduced stiffness, since the numerical sensitivities are only correct if the step size is smaller than the smoothing parameter \cite{Stewart2010}. 
This requires large optimization problems even for moderate accuracy. 
Alternatively, high accuracy event-based methods \cite{Kirches2006,Stewart1990b} are, due to external switch detecting algorithms and internal logical statements, very difficult to incorporate into direct optimal control~\cite{Kirches2006}.
In this paper, we overcome these difficulties by applying the recently introduced Finite Elements with Switch Detection (FESD) method \cite{Nurkanovic2022} to the time-freezing Filippov system. 
This method automatically detects the switching events without an external routine and delivers correct numerical sensitivities~\cite{Nurkanovic2022}.
It was demonstrated in several benchmarks \cite{Nurkanovic2022a,Nurkanovic2022b,Nurkanovic2022} that FESD outperforms standard discretization and mixed-integer approaches while providing several orders of magnitude more accurate for the same computational effort.

The time-freezing reformulation was first introduced in \cite{Nurkanovic2021} for CLS with partially elastic impacts. 
In this paper, we extend these ideas to the inelastic case. 
Time-freezing was later independently introduced by Halm and Posa \cite{Halm2021}.
They treat the inherent nonuniqueness in the simulation of CLS with multiple frictional impacts in a stochastic way. 
The multiple outcomes are achieved by randomly varying parameters of the auxiliary dynamics and by considering a set of possible simulation solutions.
In their work, time-freezing is used to obtain a differential inclusion with a bounded right-hand side, to be able to apply standard solution existence results.
In this paper, our goal is to transform an Optimal Control Problem (OCP) with CLS into an equivalent OCP subject to a Filippov system, which we can solve with high accuracy via FESD\cite{Nurkanovic2022}.

Other reformulation approaches use coordinate transformations \cite{Kim2014,Zhuravlev1978}, penalization/smoothing of the complementary \cite{Stewart2010,Tassa2012} and compliant contact models \cite[Chapter 2]{Brogliato2016}. 
The former have the advantage that they are exact, i.e., we can recover the solution of the original system.
However, they are usually limited to special cases, namely, partially elastic impacts and a single scalar constraint.
Examples are the Zhuravlev-Ivanov transformation \cite{Zhuravlev1978},\cite[Sec 1.4.3.]{Brogliato2016} and the gluing function approach \cite{Kim2014} in the hybrid systems formalism. 
\vspace{-0.2cm}
\subsection{Contributions}
\vspace{-0.2cm}
This paper extends the time-freezing reformulation from \cite{Nurkanovic2021}, which exactly transforms a CLS with a single unilateral constraint into a piecewise smooth system.
We discuss how to select auxiliary dynamics and formalize the relationship between the time-freezing system and CLS. 
Additionally, we present an extension to handle state jumps in tangential directions resulting from friction. 
We demonstrate how to apply the time-freezing system in optimal control and show that its solutions are also optimal for the initial OCP with a CLS.
Moreover, we show how to reformulate the time-freezing system into a dynamic complementarity system, which facilitates the application of the FESD method for direct optimal control.
Furthermore, we introduce time transformations and constraints to achieve equidistant control grid discretization and desired final time despite the nonsmooth clock state. 
The discretized optimal control problems result in mathematical programs with complementarity constraints. 
They are solved via a homotopy approach. 
Thereby, we solve only a few smooth Nonlinear Programs (NLP) and recover a highly accurate nonsmooth solution approximation with state jumps.
The theoretical considerations and efficacy of the proposed numerical methods are demonstrated on a challenging OCP example where the dynamic trajectory of a robot hopping over holes is computed.
All methods and examples from this paper, including a fully-automated reformulation of the CLS into a PSS, are implemented in the open-source tool NOSNOC \cite{Nurkanovic2022c,Nurkanovic2022b}.

\vspace{-0.2cm}
\subsection{Outline}
\vspace{-0.1cm}
The remainder of this paper is structured as follows. Section \ref{sec:nonsmooth_dynamics} gives an introduction to piecewise smooth systems, their embedding into Filippov's framework, and complementarity Lagrangian systems with state jumps. This is followed by Section \ref{sec:time_freezing} where the time-freezing reformulation is discussed in detail. In Section \ref{sec:frictional_impact} we extend these ideas to the frictional impact case.
Section \ref{sec:computational_considerations} discusses numerical optimal control for time-freezing systems.  Section \ref{sec:numerical_examples} provides a numerical example of a robotics optimal control problem. We conclude and list some future research directions in Section \ref{sec:conclusions}. 
\vspace{-0.2cm}
\subsection{Notation}
\vspace{-0.1cm}
Time derivatives of a function $x(t)$ w.r.t. to {the physical time} $t$ are compactly denoted by $\dot{x}(t) \coloneqq \frac{\dd x(t)}{\dd t}$,
and of a function $y(\tau)$ w.r.t. to {the numerical time} $\tau$ by ${y'(\tau) \coloneqq \frac{\dd y(\tau)}{\dd \tau}}$. 
For the left and the right limits, we use the notation ${x(\ts^+)  = \lim\limits_{t\to \ts,\ t>\ts} x(t)}$ and ${x(\ts^-)  = \lim\limits_{t\to \ts,\  t<\ts}x(t)}$, respectively.
For ease of notation, when clear from the context we drop the $t$, $\tau$, or $x$-dependencies.
All vector inequalities are to be understood element-wise. 
The complementarity conditions for two vectors  $a,b \in \R^{n}$ read as ${0\leq a \perp b\geq 0}$, where $a \perp b$ means $a^{\top}b =0$.
The matrix $I_{n} \in \R^{n\times n}$ is the identity matrix, and $\mathbf{0}_{m,n} \in \R^{m\times n}$ is a matrix whose entries are all zeros. 
The concatenation of two column vectors $a\in \R^m$, $b\in \R^n$ is denoted by $(a,b)\coloneqq[a^\top,b^\top]^\top$. 
The concatenation of several column vectors is defined analogously.
A vector with all ones is denoted by $e=(1,1,\dots,1) \in \R^n$, and its dimension is clear from the context. The closure of a set $X$ is denoted by $\overline{X}$, its boundary by $ \partial X$. 

The set-valued sign function is defined as 
\begin{align*}
	\sign(x) = \begin{cases}
		\{1\}, &x>0,\\
		[-1,1],& x=0,\\
		\{-1\},&x<0.
	\end{cases}
\end{align*}
The vector-valued version $\Sign: \R^n \rightrightarrows\R^n$ is defined as $\Sign(x) = (\sign(x_1),\ldots,\sign(x_n))$.
Key symbols and definitions used in this paper are summarized in Table~\ref{tab:symbols}.

\section{Nonsmooth differential equations} \label{sec:nonsmooth_dynamics}
\begin{table}[t!]	
	\centering
	\caption{Key symbols used throughout this paper.
	  }
	\centering
	\begin{tabular}{p{1.8cm}p{6.7cm}}	
		\hline
		\textbf{Symbol} & \textbf{Meaning and reference}\\
		\hline
		$n_x$ & dimension of the differential state $x$, Sec. \ref{sec:nonsmooth_dynamics}\\
		$n_f$ & number of regions/modes in the PSS, Eq.\eqref{eq:pws}\\
		$n_q$ & dimension of the states $q$ and $v$, Sec. \ref{sec:cls}\\
		$n_u$ & dimension of the control function, Sec. \ref{sec:cls}\\
		$n_y$ & dimension of the time-freezing state, Sec. \ref{sec:time_freezing_defintion}\\
		$\nt$ & dimension of the tangent space, Sec. \ref{sec:cls_friction}\\
		\hline
		$x$ & differential state, Sec. \ref{sec:nonsmooth_dynamics}.\\
		$\theta$ & Filippov's convex multipliers, Eq. \eqref{eq:FilippovDI_with_multiplers}\\
		$q$ & \textit{position} state, Eq. \eqref{eq:cls} \\
		$v$ & \textit{velocity} state, Eq. \eqref{eq:cls} \\
		$\vt$ & tangential \textit{velocity} at contact points, Sec.~\ref{sec:cls_slip_stick}\\
		$u$ & control function, Eq. \eqref{eq:cls} \\
		$\zn$ & Lagrange multiplier, norm. contact force, Eq.~\eqref{eq:cls} \\
		$\zt$ & Lagrange multiplier, friction force, Eq.~\eqref{eq:cls_friction} \\
		$y$ & extended state of time-freezing system, Sec. \ref{sec:time_freezing_defintion}\\
		$t$ & physical time, Sec. \ref{sec:time_freezing_defintion}\\
		$\tau $ & numerical time, Sec. \ref{sec:time_freezing_defintion}\\
		$\anormal $ & constant of the auxiliary dynamics, Prop. \ref{prop:aux_forming_ode}\\
		$\atangent $ & constant of the auxiliary dynamics, Eq. \eqref{eq:auxiliary_dynamics_tangetnial}\\
		$\lambdan$,$\lambdap$ & dual variables of the linear program, Eq. \eqref{eq:time_freezing_dcs}\\
		$\alpha$ & primal variables of the linear program, Eq.~\eqref{eq:parametric_lp}\\
		$\mu$ & coefficient of friction, Eq.~\eqref{eq:cls_friction_disipiation2}\\
		$\epst$ & relaxation parameter in friction model, Eq.~\eqref{eq:friction_slution_map_approx}\\
		$s$ & speed of time control variable, Eq. \eqref{eq:ocp_time_freezing}\\
		\hline
		$f_i(x,u)$ & modes of the PSS system, Eq. \eqref{eq:pws}\\
		$f_v(q,v,u)$ & vector field of the \textit{velocity} state, Eq.
		\eqref{eq:cls_v_dyn}\\
		$n(q)$ & normal to the CLS constraint surface, Eq. \eqref{eq:cls_v_dyn}\\
		$M(q)$ & inertia matrix, Eq.~\eqref{eq:cls_v_dyn}\\
		$f_c(q)$ & constraint function in the CLS, Eq. \eqref{eq:cls_comp}\\
		$f_{\mathrm{ODE}}(x,u)$ &  unconstrained dynamics when $f_c(q)\!>\!0$, Eq.~\eqref{eq:cls_free_flight}\\
		$D(q)$ &  Delassus' matrix/scalar, Eq.~\eqref{eq:contact_LCP_functions_D}\\
		$\varphi(x,u)$ &  determines if contact persists, Eq~\eqref{eq:contact_LCP_functions_psi}\\
		$f_{\mathrm{DAE}}(x,u)$ &  dynamics equivalent to the DAE \eqref{eq:cls_dae} (constrained dynamics when $f_c(q)=0$), Eq.~\eqref{eq:cls_index_reduced_dae}\\
		$c_i(y)$ & time-freezing PSS switching functions, Sec. \ref{sec:time_freezing_defintion}, Eq.~\eqref{eq:time_freezing_regions} and Sec. \ref{sec:time_freezing_friction}, Eq.~\eqref{eq:time_freezing_regions_friction}\\
		$f_{\mathrm{aux,n}}(y)$ & auxiliary dynamics, Def. \ref{def:auxiliary_dynamics} and Prop. \ref{prop:aux_forming_ode} \\
		$\gamma(x,u)$ &time-rescaling factor of the time-freezing sliding mode, Eq.~\eqref{eq:time_freezing_rescaling}\\
		$\Tan$ & matrix whose columns span the tangent space at contact points, Eq.~\eqref{eq:cls_friction}\\ 
		$b_j(q)$ & $j$-th column of $\Tan$, Eq.~\eqref{eq:cls_friction}\\ 
		$f_{\mathrm{aux,t}}^-(y)$ & auxiliary dynamics for tan. directions, Eq.~\eqref{eq:auxiliary_dynamics_tangetnial} \\
		$\fslip(x,u)$  & dynamics for slipping motion in contact phases, Eq.~\eqref{eq:cls_slip}\\
		$\fstick(x,u)$  & dynamics for sticking motion in contact phases, Eq.~\eqref{eq:cls_stick}\\
		$\tilde{D}(q)$  & generalization of $D(q)$, Sec.~\ref{sec:cls_slip_stick}\\
		$\tilde{\varphi}(x,u)$  & generalization of $\varphi(x,u)$, Sec.~\ref{sec:cls_slip_stick}\\
		$g_{\mathrm{F}}(\theta,\alpha)$  & expression for relating $\theta$ and $\alpha$, Eq.~\eqref{eq:time_freezing_dcs_lifting}\\
		$\Psi(x(T))$  & terminal cost, Eq.~\eqref{eq:ocp_cls}\\
		$g(x,u)$  & path and terminal constraints, Eq.~\eqref{eq:ocp_cls}\\
		$r(x)$  & terminal constraints, Eq.~\eqref{eq:ocp_cls}\\
		\hline
		$R_i$ & regions of the PSS, Eq.~\eqref{eq:pws}\\
		$\mathcal{I}$ & index set for PSS modes, Eq.~\eqref{eq:pws}\\
		$F_{\mathrm{F}}(x,u)$  & Filippov set, Eq.~\eqref{eq:FilippovDI_with_multiplers}\\
		$\Sigma$ & switching surface, Sec.~\ref{sec:time_freezing_defintion}\\
		$F_{\mathrm{TF}}(y,u)$  & Filippov set for the time-freezing system, Def.~\ref{def:time_freezing_system} and \ref{def:time_feezing_friction}\\
		$Q$  & region where all auxiliary dynamics are defined, Eq.~\eqref{eq:time_freezing_regions_friction}\\
		\hline
	\end{tabular}
	\vspace{-0.20cm}
	\label{tab:symbols}
\end{table}

In this section, we define piecewise smooth systems (PSS), Filippov's notion of solutions for the PSS \cite{Filippov1988}, and complementarity Lagrangian systems (CLS).
\subsection{Piecewise smooth systems}\label{sec:pss_and_filippov}
We regard piecewise smooth systems of the following form:
\begin{align} \label{eq:pws}
	\dot{x} & =f_i(x,u),\ \text{if}\;  x\! \in\! R_i \subset \R^{n_x},\ i \in\! \mathcal{I}\! \coloneqq\! \{ 1,\ldots,n_f \},
\end{align}
with regions $R_i\subset \R^{n_x}$ and associated dynamics $f_i(\cdot)$, which are at least twice continuously differentiable functions on an open neighborhood of $\overline{R}_i$. The control function $u \in \R^{n_u}$ is assumed to be given and can be obtained, e.g., by solving an optimal control problem.
The right-hand side of \eqref{eq:pws} is in general discontinuous in $x$.
We assume that the sets $R_i$ are disjoint, nonempty, connected, and open. They have piecewise smooth boundaries $\partial R_i$. Moreover, it is assumed that $\overline{\bigcup\limits_{i\in \mathcal{I}} R_i} = \R^n$ and that $\R^n \setminus \bigcup\limits_{i\in\mathcal{I}} R_i$ is a set of measure zero. 

The ODE \eqref{eq:pws} is not properly defined on the boundaries $\partial R_i$. To have a meaningful solution concept for the PSS \eqref{eq:pws} we regard its Filippov extension~\cite{Filippov1988}. The ODE \eqref{eq:pws} is replaced by a differential inclusion whose r.h.s. is a convex and bounded set. Due to the assumed structure of the sets $R_i$, if $\dot{x}$ exists, functions $\theta_i(\cdot)$ that serve as convex multipliers, can be introduced, and the Filippov differential inclusion for \eqref{eq:pws} reads as \cite{Filippov1988,Stewart1990b}:
\begin{align}\label{eq:FilippovDI_with_multiplers}
	\begin{split}
		\dot{x}  \in    F_{\mathrm{F}}(x,u) \coloneqq \Big\{&\sum_{i\in \mathcal{I}}
		f_i(x,u) \, \theta_i  \mid \sum_{i\in \mathcal{I}}\theta_i = 1,
		\ \theta_i \geq 0,\\
		&0= \theta_i \  \mathrm{if} \;  x \notin \overline{R_i}, 
		\forall  i  \in \mathcal{I} \Big\}.
	\end{split}
\end{align}
Note that in the interior of the regions $R_i$ the \textit{Filippov set} $F_\mathrm{F}(x,u)$ is equal to $\{f_i(x,u)\}$ and on the boundary between regions it is a convex combination of the neighboring vector fields.
Sufficient conditions for the existence and uniqueness of solutions are given in \cite{Filippov1988}. 
We assume that the boundaries of the regions $\partial R_i$ are defined by the zero-level sets of, at least twice continuously differentiable, scalar functions $c_i(x) = 0$. These functions are called \textit{switching functions}.
The evolution of $x(\cdot)$ on region boundaries $\partial R_i$, are called \textit{sliding modes}. The dynamics of sliding modes are implicitly defined by differential
algebraic equations, since the corresponding constraints $c_i(x)= 0$ must hold~\cite{Filippov1988}.

\subsection{Complementarity Lagrangian systems}\label{sec:cls}
This paper regards complementarity Lagrangian systems (CLS) with a single unilateral constraint. 
A CLS with inelastic impacts reads as  
	\begin{subequations}\label{eq:cls}
		\begin{align}
			\dot{q} &= v, \label{eq:cls_q_dyn} \\
			\dot{v} &= f_\mathrm{v}(q,v,u)  + M(q)^{-1} n(q) \zn  \label{eq:cls_v_dyn},\\
			0 &\leq \zn \perp f_c(q) \geq 0 \label{eq:cls_comp}, 
			\\
			\begin{split}
				0&=n(q(\ts))^\top v(\ts^+),\\
				&\mathrm{if} \ f_c(q(\ts)) = 0\ \mathrm{and}\ n(q(\ts))^\top v(\ts^-)<0,
			\end{split} \label{eq:cls_state_jump_law}
		\end{align} 
	\end{subequations}
where $x \coloneqq (q,v)$ and the states $q\in \R^{n_q}$ and $v\in \R^{n_q}$ correspond to the position and velocity of a rigid body, respectively. 
The variable $\zn\in \R$ is the Lagrange multiplier and is physically interpreted as the normal contact force.
The function $f_c(q)\in \R$ is the signed distance between a contact point of a rigid body and an obstacle or another rigid body, and $n(q) \coloneqq \nabla_q f_c(q)$.
The matrix $M(q)$ is the inertia matrix and the function $f_\mathrm{v}(q,v,u)$ corresponds to the total \textit{acceleration} of the rigid body, i.e., it collects all internal and external forces (except the contact forces), multiplied by the inverse of $M(q)$.
The functions $f_\mathrm{v}: \R^{n_q} \times \R^{n_q} \times \R^{n_u} \to \R^{n_q}$, ${M(q):\R^{n_q}  \to \R^{n_q \times n_q}}$, $f_c: \R^{n_q} \to \R$ are assumed to be at least twice continuously differentiable and the matrix $M(q)$ is assumed to be symmetric positive definite. 
The complementarity condition \eqref{eq:cls_comp} states that: either the system is in contact and there is a reaction force ($f_c(q)=0,\; \zn\geq0$) or there is no contact and no reaction force ($f_c(q)>0,\; \zn=0$). 
When the body makes contact, the negative normal velocity $n(q)^\top v$ must jump to zero. 
This is modeled via Eq. \eqref{eq:cls_state_jump_law}.

The function $u\in \R^{n_u}$ in Eq.~\eqref{eq:cls_v_dyn} is the control function. 
In this paper, we aim to find the control function $u(t)$ by solving the Optimal Control Problem (OCP):
\begin{subequations} \label{eq:ocp_cls1}
	\begin{align}
		\min_{x(\cdot),\zn(\cdot),u(\cdot),} \quad &   \Psi(x(\Tctrl)) \label{eq:ocp_cls1_objective}\\
		\textrm{s.t.} \quad  &x(0) = \bar{x}_0 \label{eq:ocp_cls1_iv},\\
		&\textrm{Eq.} \eqref{eq:cls},\; t \in [0,\Tctrl] \label{eq:ocp_cls1_dynamics}\\
		&0\leq g(x(t),u(t)),\; t \in [0,\Tctrl] \label{eq:ocp_cls1_path},\\
		&0\leq r(x(\Tctrl)),\;  \label{eq:ocp_cls1_temrinal}
	\end{align}
\end{subequations}
where $\Psi: \R^{n_x} \to \R$ is the terminal cost of the OCP, $\bar{x}_0$ is a given initial value. 
The functions $g: \R^{n_x} \times \R^{n_u} \to \R^{n_g}$ and $r: \R^{n_x} \to \R^{n_r}$ are the path and terminal constraints, respectively. 
The CLS dynamics in Eq. \eqref{eq:ocp_cls1_dynamics} (resp. Eq. \eqref{eq:cls}) make this OCP nonsmooth and nonconvex. 

The nonsmooth dynamics and discontinuous velocity state make this OCP difficult to solve numerically. 
In the subsequent sections, we develop the time-freezing reformulation such that we can reformulate the problem \eqref{eq:ocp_cls1} into an OCP subject to a Filippov system, for which more efficient numerical methods are available~\cite{Nurkanovic2022}.
In this paper, we focus on the case of a single unilateral constraint $f_c(q)$. 
Extensions for multiple and simultaneous impacts will be studied in future work. 
The extension depends on the chosen impact model \cite{Brogliato2016}.

Next, we introduce a guiding example on which we will illustrate the main ideas behind the time-freezing reformulation throughout the paper.
\begin{example}[Guiding example]\label{ex:guiding_example}
Consider a frictionless point mass in two dimensions above a horizontal table. The mass is $m=1$ $\mathrm{kg}$ and  $g = 9.81\, {\mathrm{m}}/{\mathrm{s}^2}$ is the gravitational acceleration. Denote by $q \coloneqq (q_1,q_2)$ and $v \coloneqq (v_1,v_2)$ its position and velocity, respectively, and let $\zn$ be the normal contact force. The dynamics are given by the CLS:
	\begin{subequations}\label{eq:example_cls}
		\begin{align}
			\dot{q} &= v,\\
			m\dot{v} &=  \begin{bmatrix}	0 \\ -mg\end{bmatrix}
			 + \begin{bmatrix}	0 \\ 1 \end{bmatrix} \zn
			 +  \begin{bmatrix}	u_1 \\ u_2	\end{bmatrix} , \\ 
		0 &\leq \zn \perp q_2 \geq 0, \label{eq:example_complement}\\
		v_2(\ts^+)& \!=\! 0,  \ \mathrm{if} \ q_2(\ts) \!=\! 0\
				\mathrm{and}\ v_2(\ts^-)\!<\!0.\label{eq:example_state_jump_law}
		\end{align}
	\end{subequations}
where $u = (u_1,u_2)\in \R^2$ is an externally chosen thrust force, which shall be found, e.g., by solving an optimal control problem.	
\end{example}
\section{The time-freezing reformulation} \label{sec:time_freezing}

This section develops the time-freezing reformulation that enables one to transform the CLS \eqref{eq:cls} with inelastic impacts into a PSS of the form of \eqref{eq:pws}. For simplicity and ease of exposition, we first focus on the case without friction. Extensions with frictional impacts are given in Section \ref{sec:frictional_impact}. We start by investigating different possible modes of the CLS  \eqref{eq:cls}. 
Afterward, the time-freezing reformulation with its needed ingredients is introduced. 
The section finishes by formally relating the two regarded systems.

\subsection{The different modes of the CLS}\label{sec:cls_modes}
For the CLS \eqref{eq:cls} we can distinguish two modes of operation:

(i) the system is not in contact (unconstrained case, free flight), i.e., $f_c(q)>0$ which implies $\zn = 0$,
(ii) the system is in contact, i.e., $f_c(q)=0$ and $\zn\geq 0$. In the first case, the system evolves according to the ODE:

	\begin{align}\label{eq:cls_free_flight}	 
		\dot{q} = v,\;\dot{v} &= f_{\mathrm{v}}(q,v,u). 
	\end{align}
We write this ODE compactly as $ \dot{x} = {f}_{\mathrm{ODE}}(x,u) \coloneqq (v,f_{\mathrm{v}}(q,v,u))$.

We call an active-set change from \sloppy$\zn(\ts^-) = 0,\ f_c(q(\ts^-)) \geq 0$ to $\zn(\ts^+) \geq 0,\; f_c(q(\ts^+))  = 0$, which triggers a state jump, an \textit{impact}.

After an impact, it holds that $0=n(q(\ts))^\top v(\ts^+)$. Subsequently, the system evolves according to a differential algebraic equation (DAE) of index~3:
\begin{subequations}\label{eq:cls_dae}
	\begin{align}	
		\dot{q} & = v \label{eq:cls_dae_differential2},\\
		\dot{v} &= f_{\mathrm{v}}(q,v,u)+  M(q)^{-1} n(q) \zn, \label{eq:cls_dae_differential1}\\
		0& = f_c(q) \label{eq:cls_dae_algebraic}.
	\end{align}
\end{subequations}
Note that $f_c(q(t))$ needs to be differentiated twice w.r.t. to time until $\zn(t)$ appears explicitly. The next question to be answered is: will the system stay in contact (dynamics defined by \eqref{eq:cls_dae} with $ f_c(q)=0$) or will the \textit{contact break} (dynamics defined by \eqref{eq:cls_free_flight} with $f_c(q)>0$)?  The answer can be found by looking at the \textit{contact Linear Complementarity Problem} (LCP) \cite[Section 5.1.2]{Brogliato2016}. 
Under our standing assumptions, during contact on some time interval $[t_1,t_2]$ the \textit{consistent initialization} conditions hold
	\begin{align}\label{eq:consistent_init}
		0& = f_c(q(t)),\	0 = \dt f_c(q(t)) = \nabla_q f_c(q(t))^\top v(t). 
	\end{align} 
Consequently, $\zn(t) \geq 0, t \in [t_1,t_2]$. Due to the continuity of $q(t)$, $f_c(q(t))$ and $\dt f_c(q(t))$, for contact breaking (i.e., $f_c(q)$ becomes strictly positive) it is required that $\frac{\dd^2}{\dd t^2} f_c(q(t)) \geq 0$ for $t \in \left[t_2,t_2+\hat{\epsilon}\right)$, for some $\hat{\epsilon}>0$. Therefore, from \eqref{eq:cls_comp} we deduce that
\begin{align}\label{eq:contact_mLCP}
	0 & \leq \frac{\dd^2}{\dd t^2} f_c(q(t)) \perp \zn(t) \geq 0,\; t\in [t_1,t_2].
\end{align}
Then, by computing $\frac{\dd^2}{\dd t^2} f_c(q(t)) $ and using the r.h.s. of \eqref{eq:cls_dae_differential1}, we obtain the contact LCP in $\zn(t)$: 
\begin{align}\label{eq:contact_LCP}
	0 & \leq D(q) \zn + \varphi(x,u) \perp \zn \geq 0,
\end{align}
with
\begin{subequations}\label{eq:contact_LCP_functions}
	\begin{align}
		D(q) &= \nabla_q f_c(q)^\top M(q)^{-1} \nabla_q f_c(q), \label{eq:contact_LCP_functions_D}\\
		\varphi(x,u) &= \nabla_q f_c(q)^\top f_{\mathrm{v}}(q,v,u) + \nabla_q(\nabla_q f_c(q)^\top v)^\top v \label{eq:contact_LCP_functions_psi},
	\end{align}
\end{subequations}
where $D(q)>0$~\cite{Brogliato2016}.

The solution map of the LCP \eqref{eq:contact_LCP} is given by
\begin{align}\label{eq:contact_LCP_solution_map}
		\zn &= \max(0,-D(q)^{-1}\varphi(x,u) ).
\end{align}
From the last equation we deduce that contact breaking or sticking depends on the sign of the function  $\varphi(x,u)$.

In the case of $\varphi(x,u)\leq 0$ from Eq. \eqref{eq:contact_LCP} and \eqref{eq:contact_LCP_solution_map} it follows that $\zn(t) \geq 0$ and $\frac{\dd^2}{\dd t^2} f_c(q(t)) = 0$. Therefore, we have a persistent contact and the system evolves according to the DAE \eqref{eq:cls_dae}. Using index reduction and the solution map \eqref{eq:contact_LCP_solution_map} we can derive an ODE that is equivalent to the DAE \eqref{eq:cls_dae}:
\begin{subequations}\label{eq:cls_index_reduced_dae}
\begin{align}
\dot{q} &= v,\\
\dot{v} &=  f_{\mathrm{v}}(q,v,u) - ~ M(q)^{-1} n(q) D(q)^{-1}\varphi(x,u).
\end{align}
\end{subequations}
We compactly denote this ODE by $\dot{x} = f_{\mathrm{DAE}}(x,u)$.

In the second case, $\varphi(x,u) > 0$ implies $\zn(t) = 0$ and $\frac{\dd^2}{\dd t^2} f_c(q(t))>0$, therefore the contact breaks and the system evolves according to the ODE \eqref{eq:cls_free_flight}.

To summarize, if the system switches from the ODE mode in \eqref{eq:cls_free_flight} to the DAE mode in \eqref{eq:cls_index_reduced_dae}, a state jump must occur, except if the active-set changes happen with $n(q(\ts^-)^\top v(\ts^-)=0$. Now the system evolves on the boundary of the feasible set with $f_c(q) = 0$ according to the DAE \eqref{eq:cls_dae}, or equivalently according to the ODE defined by \eqref{eq:cls_index_reduced_dae}.
On the other hand, if we switch from DAE to ODE mode, we have a continuous transition without state jumps, i.e., {contact breaking} occurs.
\subsection{Main ideas and auxiliary dynamics}\label{sec:time_freezing_defintion}
The arguments above reveal that the CLS \eqref{eq:cls} switches between an ODE and a DAE of index 3. 
This already bears similarity to a PSS, but the main obstacle to completing this transition are the state jumps. Note that large parts of the state space, namely $f_c(q)<0$, are prohibited for the solution trajectories of the CLS.

The time-freezing reformulation is based on the following two main ideas~\cite{Nurkanovic2021}. First, we relax the constraint and allow $f_c(q)<0$. We define an \textit{auxiliary dynamical system} in this \textit{infeasible} region whose trajectory endpoints satisfy the state jump law \eqref{eq:cls_state_jump_law} on some finite time interval.  Second, we introduce a \textit{clock state} $t(\tau)$ that stops counting (i.e., $t'(\tau) = 0$), when the auxiliary ODE is active (for $f_c(q)< 0$). By taking the pieces of the trajectory when the clock state was active, one can recover the solution of the original system with discontinuous trajectories. Note that the time-freezing system has no discontinuities in its solution, but only in its r.h.s..
The extended state of the time-freezing system reads as $y  \coloneqq (x,t) \in \R^{n_y}$, $n_y = n_x+1$. 
The time of the time-freezing system $\tau$ is called \textit{numerical time.} The intervals with $t'(\tau) > 0$ are referred to as \textit{physical time} and those with $t'(\tau)=0$ as \textit{virtual time}.
The properties of the auxiliary dynamics are summarized in the following definition. 
Later a constructive way to find such systems is given.

\begin{definition}[Auxiliary dynamics]\label{def:auxiliary_dynamics}
An auxiliary dynamical system ${y}'(\tau) = f_{\mathrm{aux,n}}({y}(\tau))$ satisfies for every initial value ${y}(\taus) = (\qs,\vs,t_s)$, with $f_c(\qs)=0$ and $n(\qs)^\top \vs <0$,  for every well-defined and finite time interval $(\taus,\taur)$, with the length $\tau_{\mathrm{jump}} = \taur -\taus$, the following properties:
(i) $f_c(q(\tau)) \leq 0, t'(\tau) =0 \ \forall \tau \in (\taus,\taur)$,
(ii) $n({q}(\taur))^\top {v}(\taur)=0$, and
(iii) $f_c({q}(\taur)) = 0$.
\end{definition}

To construct the {time-freezing system} we take several steps.
First, observe that the post-impact velocity \eqref{eq:cls_state_jump_law} is equal to the total time derivative of the constraint $f_c(q)=0$, i.e., $\dt f_c(q) = n(q)^\top v = 0$. We choose these functions as switching functions, i.e., $c_1(y) = f_c(q)$ and $c_2(y) = n(q)^\top v$, and define the following regions:
	\begin{align}
		\begin{split}\label{eq:time_freezing_regions}
			&R_1 \!= \! \{y \! \in \! \R^{n_y} \!  \mid \!  c_1(y)\! >\! 0\}
	\! \cup\!  \{y \! \in \! \R^{n_y} \! \mid \!  c_1(y)\! < \! 0, c_2(y)\! >\! 0 \},\\
		&R_2 \!= \!\{y \in \R^{n_y} \mid c_1(y)<0, c_2(y)<0\}.
	\end{split}
	\end{align}
Second, we associate with the region $R_1$ the \textit{unconstrained dynamics} $y' = (f_{\mathrm{ODE}}(x,u),1)$ and with 
$R_2$ the auxiliary dynamics from Definition \ref{def:auxiliary_dynamics}. 
Note that the {unconstrained dynamics} is Eq. \eqref{eq:cls_free_flight} augmented by the clock state dynamics $t' = 1$.
Moreover, the control functions $u(\cdot)$ does not influence the dynamics for $y \in R_2$, i.e., whenever the time is time frozen the control can take any value.

We formally define the time-freezing system and its Filippov extension in the next definition.
\begin{figure}[t]
	\centering
	\centering
	\vspace{0.05cm}
	{\includegraphics[scale=0.60]{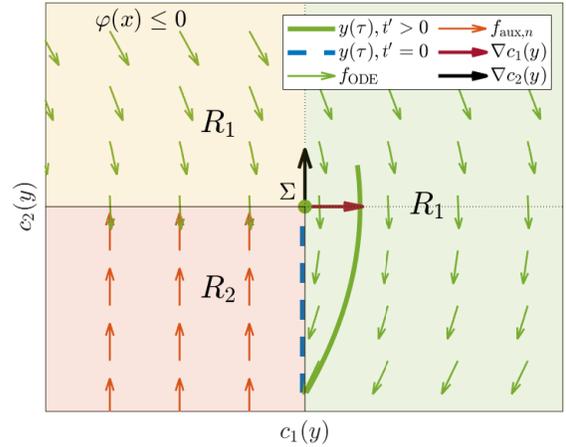}}
	\vspace{-0.6cm}
	\caption{Illustration of a phase plot of the time-freezing system from Definition  \ref{def:time_freezing_system}. The red and yellow shaded areas are infeasible for the CLS \eqref{eq:cls}. The trajectories of the auxiliary dynamics (the blue dashed line) flow in the red-shaded area.
	}
	\label{fig:time_freezing_2d}
\end{figure}
\begin{definition}[Time-freezing system]\label{def:time_freezing_system}
	Let $\tau \in \R$ be the numerical time and $y(\tau) \coloneqq ({x}(\tau),t(\tau)) \in \R^{n_y}$ the differential states and $u(\tau)\in\R^{n_u}$ a given control function. 
	The time-freezing PSS is a PSS of the form of Eq. \eqref{eq:pws}, defined by the regions $R_1$ and $R_2$ in \eqref{eq:time_freezing_regions} with $f_1(y,u) = (f_{\mathrm{ODE}}(x,u),1)$ and $f_2(y) = f_{\mathrm{aux,n}}(y)$.
	The corresponding Filippov system, which we call the time-freezing system, is defined as
	\begin{align}\label{eq:time_freezing_system}
		\begin{split}
		\!\! y' \!\in  \!F_{\mathrm{TF}}(y,u) \! \coloneqq \! \Big\{&\theta_1 {f}_1(y,u) + \theta_2 f_2(y) \mid 
		e^\top \theta \! =\! 1, \theta \!\geq  \!0 \\
		&\theta_i = 0, \; \mathrm{if} \; y \notin \overline{R_i}, i =1,2
		 \Big\}\!,
	 		\end{split}
\end{align}
with $\theta = (\theta_1,\theta_2)$. It is assumed that appropriate dynamics $f_{\mathrm{aux,n}}({y})$ exist.
\end{definition}

The phase plot of the time-freezing system is depicted in Fig. \ref{fig:time_freezing_2d}.  Note that the region $R_1$ consists of $R_1^a = \{y \in \R^{n_y} \mid  f_c(q)>0\}$ (green area), that corresponds to the feasible set of the unconstrained dynamics \eqref{eq:cls_free_flight} and the set $R_{1}^b \coloneqq \{y \in \R^{n_y} \mid c_1(y)<0, c_2(y)>0 \}$ (yellow area). The solution trajectories never flow in $R_1^b$ and the system should not be initialized in there. However, as we show later, it is crucial for sliding modes and contact breaking.
Region $R_2$ (red shaded area) contains the auxiliary dynamics that mimic the state jump.

To make further use of Definition \ref{def:time_freezing_system} we must specify how to select an appropriate auxiliary ODE. The next proposition provides a constructive way of selecting the auxiliary ODE from Definition \ref{def:auxiliary_dynamics} for any smooth scalar constraint $f_c(q) = 0$.
\begin{proposition}[Auxiliary dynamics]\label{prop:aux_forming_ode}
	Suppose that  $y(\taus)= (\qs,\vs,\ts)$ is given  such that $f_c(\qs)=0$ and $n(\qs)^\top \vs \leq 0$ holds. Then the ODE given by
	\begin{align} \label{eq:aux_dyn_example}
		{y}' = f_{\mathrm{aux,n}}(y) & 	\coloneqq \begin{bmatrix}
			\mathbf{0}_{n_q , 1} \\   M({q})^{-1} n({q}) \anormal \\0 
		\end{bmatrix}
	\end{align}
	with $\anormal>0$  is an auxiliary dynamical system from Definition \ref{def:auxiliary_dynamics} with~${\taujump= -\frac{n(\qs)^\top \vs}{D(\qs)\anormal}}$.
\end{proposition}
\textbf{PROOF.}
According to \eqref{eq:aux_dyn_example} we have ${q}'(\tau)=  \mathbf{0}_{n_q,1},\ \forall \tau \geq \taus$, which implies ${q}(\tau)  = \qs$ and $f_c({q}(\tau))= 0,\ \forall \tau \geq \taus$.
This means also that $M({q}(\tau)) = M(\qs),\ \forall \tau \geq \taus$. 
Second, regard the dynamics of ${v}' = M({q})^{-1} n({q}) \anormal = M(\qs)^{-1} n(\qs) \anormal$ and rewrite this equation in integral form. By multiplying it from the left by $n(\qs)^\top$ we obtain:
\begin{align*}
	n(\qs)^\top {v}(\tau) = n(\qs)^\top \vs \! +\! n(\qs)^\top  M(\qs)^{\!-\!1} n(\qs) \anormal (\tau-\taus).
\end{align*} 
Since the first term on the r.h.s. is negative and the second strictly positive, we deduce that	$n({q}(\taur))^\top {v}(\taur) = 0$  and $f_c({q}(\taur))=0$ with $ \taujump = \taur-\taus = -\frac{n(\qs)^\top \vs }{ D(\qs)\anormal}$. Hence, all conditions from Definition \ref{def:auxiliary_dynamics} are satisfied and the proof is complete.
\qed

Next, we discuss which mode of the time-freezing system matches the persistent contact dynamics of the CLS~\eqref{eq:cls}. We observe that the set $\Sigma \coloneqq  \{ y \mid c_1(y) = f_c(q) = 0, c_2(y) = n^\top v = 0 \}$ is defined by the same equations as the consistent initialization conditions \eqref{eq:consistent_init} (but now in $\R^{n_y}$ instead of $\R^{n_x}$ due to the clock state). A sliding mode of the time-freezing PSS evolves on $\Sigma$ just as the solution of the persistent contact DAE \eqref{eq:cls_dae}. Therefore, its dynamics should match the dynamics of the DAE \eqref{eq:cls_dae}. We detail in the next two subsections that this is indeed the case.

To summarize, for $f_c(q) >0$ (which is a subset of $R_1$) the time-freezing system and CLS have the same dynamics. In $R_2$ the auxiliary dynamics mimics the state jump and the sliding mode $y\in \Sigma$ should match the dynamics of the CLS in contact mode. To illustrate the developments so far, we derive a time-freezing PSS for our guiding example.
\begin{example}(Guiding example as time-freezing PSS)\label{ex:guiding_example_pss}
The state space is $y = (q,v,t) \in \R^5$ and the switching functions read as $c_1(y) = y_2$ and $c_2(y) = v_2$. The two PSS regions are $R_1 = \{y\in\R^5 \mid q_2 >0 \} \cup \{y\in\R^5 \mid q_2 <0, v_2>0 \}$ and
$R_2 = \{y\in\R^5 \mid q_2 <0, v_2<0 \}$, cf. Fig.~\ref{fig:time_freezing_2d}.
The PSS dynamics are given by $f_1 = (v_1,v_2,u_1,-g+u_2,1)$. 
The constraint normal is $n=(0,1)$, so that by applying Eq. \eqref{eq:aux_dyn_example} we find that $f_2 = (0,0,0,\anormal,0)$. 
\end{example}
\subsection{Persistent contact and sliding mode}\label{sec:presistent_contact}
Depending on the sign of the function $\varphi(x,u)$, a solution $y$ initialized at $\Sigma$ should either stay at $\Sigma$ (sliding mode, persistent contact) or leave it (contact breaking).
In this subsection, we study the case when a solution of the time-freezing system satisfies the conditions $y(\tau) \in  \Sigma$  and $\varphi(x(\tau),u(\tau)) \leq 0$ for some $\tau \in [\tau_1,\tau_2]$ (persistent contact).
During contact, the CLS system satisfies the consistent initialization \eqref{eq:consistent_init} which corresponds to $\Sigma$ without the clock state.
It is desired that under these conditions $y(\tau)$ stays on $\Sigma$ and that the corresponding sliding mode dynamics match the DAE dynamics \eqref{eq:cls_index_reduced_dae}.
For a solution to stay in the sliding mode, the surface $\Sigma$ must be \textit{stable}, i.e., all neighboring vector fields point toward $\Sigma$. Since $\nabla c_1(y)^\top f_1(y,u) = 0$, $\nabla c_1(y)^\top f_2(y) = 0$ and $\nabla c_2(y)^\top f_1(y,u) = \varphi(x,u) \leq 0$, $\nabla c_2(y)^\top f_2(y) = D(q)\anormal >0$, we see that this is indeed the case, cf. Fig.~\ref{fig:time_freezing_2d}.
We show next that the sliding mode of the time-freezing system is unique and that it matches the dynamics of the DAE of index 3 after the state jump, as required.
\begin{theorem}[Unique sliding mode]\label{th:presistent_contact}
Regard the time-freezing system from Definition \ref{def:time_freezing_system} with the auxiliary dynamics from Proposition~\ref{prop:aux_forming_ode}. Let $y(\tau)$ be a solution of this system with $y(0) \in \Sigma$ and $\tau \in [0,\tauf]$. Suppose that $\varphi(x(\tau),u(\tau)) \leq 0$ for all $\tau \in [0,\tauf]$ (persistent contact), then the following statements are true:
\vspace{-0.35cm}
	\begin{enumerate}[(i)]
		\item the convex multipliers $\theta_1,\theta_2 \geq 0$ in Eq. \eqref{eq:time_freezing_system} are unique,
		\item the dynamics of the sliding mode are given by ${y'=\gamma(x,u)(f_{\mathrm{DAE}}(x,u),1)}$, where $\gamma(x,u)\in \left(0,1\right]$ is a time-rescaling factor given by
		\begin{align}\label{eq:time_freezing_rescaling}
			\gamma(x,u) \coloneqq	\frac{D(q)\anormal}{D(q)\anormal-\varphi(x,u)}.
		\end{align}
	\end{enumerate}
\end{theorem}
\textbf{PROOF.} 
To compute the convex multipliers $\theta_1$ and $\theta_2$ we use Definition \ref{def:time_freezing_system} and the fact that $y \in \Sigma$. This results in the conditions:
\begin{align*}
c_1(y) = 0,\; c_2(y) = 0,\; \theta_1 + \theta_2 = 1.
\end{align*}
Since $\theta$ does not explicitly appear in the first two conditions, we differentiate them w.r.t. to $\tau$ and use Eq. \eqref{eq:time_freezing_system}. 
We have two unknowns and three conditions, hence the system is over-determined. However, we have by assumption that $c_2(y) = \nabla_q c_1(y)^\top  v = 0$ and by direct evaluation, we conclude that $\frac{\dd c_1(y)}{\dd \tau} = 0$ is satisfied for every $\theta_1$ and $\theta_2$.
The conditions that are left are $\frac{\dd c_2(y)}{\dd \tau} = 0$ and  $\theta_1 + \theta_2 =1$. 
Since $\nabla_v c_2(y) = \nabla_v (\nabla_q c_1(y)^\top v) = \nabla_q c_1(y)$ and $\frac{\partial c_2(y)}{\partial t} =0$ we obtain from $\frac{\dd c_2(y)}{\dd \tau} =0$ the following equations
\begin{align*}
	\begin{split}
	0 &= \frac{\dd c_2(y)}{\dd \tau} = \theta_1 
	\begin{bmatrix}
		\nabla_q c_2(y)^\top \nabla_v c_2(y)^\top
	\end{bmatrix}
	\begin{bmatrix}
		v \\ f_v(q,v,u)
	\end{bmatrix}\\
	&+\theta_2
	\begin{bmatrix}
		\nabla_q c_2(y)^\top \nabla_v c_2(y)^\top
	\end{bmatrix}
	\begin{bmatrix}
		0 \\ M(q)^{-1}\nabla_q c_1(y) \anormal
	\end{bmatrix},
	\end{split}\\
	\begin{split}
	0 &= \theta_1 \underbrace{\nabla_q (\nabla_q f_c(q)^\top v )^\top v+ \nabla_q f_c(q)^\top f_v(q,v,u) }_{= \varphi(x,u) <0} \\
	&+\theta_2 \underbrace{\nabla_q f_c(q)^\top M(q)^{-1}\nabla_q f_c(q)}_{= D(q)> 0}\anormal. 
	\end{split}
\end{align*}
Thus we obtain a system linear in $\theta$
\begin{align*}
\begin{bmatrix}
	\varphi(x,u) & D(q)\anormal \\
	1 & 1
\end{bmatrix}
\begin{bmatrix}
\theta_1\\ 
\theta_2
\end{bmatrix} &=
\begin{bmatrix}
	0\\ 1
\end{bmatrix},
\end{align*}
 and by solving it we have that $\theta_1 = \frac{D(q)\anormal}{D(q)\anormal-\varphi(x,u)} = \gamma(x,u)$ and $\theta_2 =  \frac{-\varphi(x,u)}{D(q)\anormal-\varphi(x,u)}$. Since $D(q)\anormal-\varphi(x,u)> 0$ and $\varphi(x,u)\leq 0$, we have always unique $\theta_1, \theta_2 \geq 0$. This completes the first part of the proof.
 
 For the second part, we evaluate
 \begin{align*}
y' &= \theta_1 f_1(y,u) + \theta_2 f_2(y)= \gamma(x,u)\begin{bmatrix}
	v \\ f_v(q,v,u) \\ 1 
\end{bmatrix}\\
&+\frac{-\varphi(x,u)}{D(q)\anormal-\varphi(x,u)}\!
 \begin{bmatrix}
\mathbf{0}_{n_q,1} \\ M(q)^{-1}\nabla f_c(q) \anormal \\ 0 
\end{bmatrix}\! \cdot\! D(q)^{-1} D(q). \\
\end{align*}
In the second term we use that $\zn = -D(q) \varphi(x,u)$ (cf. Eq. \eqref{eq:contact_LCP_solution_map}) and the expression for $\gamma(x,u)$ in Eq. \eqref{eq:time_freezing_rescaling}. By comparing the last expression to Eq. \eqref{eq:cls_index_reduced_dae} we obtain $y' = \gamma(x,u)(f_{\mathrm{DAE}}(x,u), 1)$.
This completes the proof. \qed

\begin{figure}[t]
	\centering
	\centering
	\vspace{-0.2cm}
	{\includegraphics[scale=0.60]{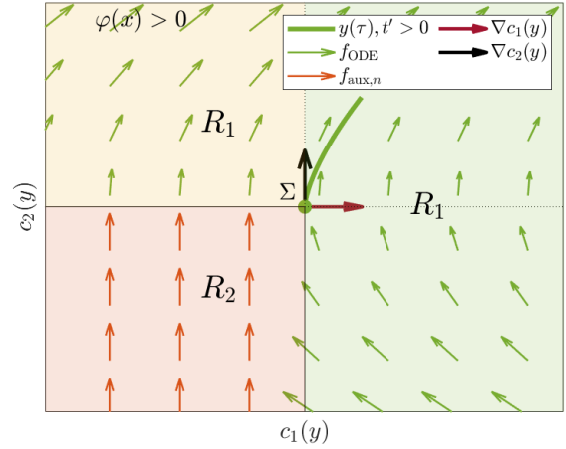}}
	\vspace{-0.7cm}
	\caption{Illustration of a phase plot of the time-freezing system from Definition  \ref{def:time_freezing_system} with $\varphi(x,u)> 0$. Compared to Fig.~\ref{fig:time_freezing_2d}, the vector field in $R_1$ is changed and $\Sigma$ is not stable anymore, thus leaving ${\Sigma}$ into $R_1$ is possible.}
	\label{fig:time_freezing_2d_exit}
\end{figure}

This theorem shows that the sliding mode of the time-freezing system on $\Sigma$ is unique and equal to the dynamics of the CLS in the persistent contact mode given by Eq. \eqref{eq:cls_index_reduced_dae}, but slowed down by the factor $\gamma(x,u)$.
Note that for larger values of $\anormal$ the factor $\gamma(x,u)$ comes closer to one, which reduces the slow-down, cf. Example \ref{ex:guiding_example_simulation}.

However, in the function $t(\tau) \to x(t(\tau)))$ the solution with a speed of time of one is recovered.
We briefly discuss the intuition behind the time slow down by $\gamma(x,u)$. 
To achieve the sliding mode on $\Sigma$, the vector fields from $R_1$ and $R_2$ (the yellow part in Fig.~\ref{fig:time_freezing_2d}) must push toward $\Sigma$. The resulting dynamics is a convex combination of the two vector fields, and since the speed of time in $R_1$ is one and in $R_2$ zero we obtain a slow down equal to $\theta_1 = \gamma(x,u)$. 
Moreover, the vector field $f_{\mathrm{ODE}}(\cdot)$ in the yellow area "stops" the trajectory coming from $R_2$ and thus enables the sliding mode. 
This shows the significance of having the vector field $(f_{\mathrm{ODE}}(\cdot),1)$ in the yellow area $R_1^b$, even though by construction the solution never flows there.

\subsection{Contact breaking}\label{sec:contact_breaking}
It is left to study the case when $y(\tau_1) \in \Sigma$ but $\varphi(x(\tau),u(\tau)) >0$ for $\tau \in [\tau_1,\tau_2]$.
In the CLS for $\varphi(x,u)>0$ the contact breaks, cf. Sec \ref{sec:cls_modes}. In the time-freezing system, we expect the trajectory to leave the sliding mode from $\Sigma$. Since $\nabla c_1(y)^\top f_1(y,u) = 0$, $\nabla c_1(y)^\top f_2(y) = 0$ and $\nabla c_2(y)^\top f_1(y,u) = \varphi(x,u) > 0$, $\nabla c_2(y)^\top f_2(y) = D(q)\anormal >0$, the surface $\Sigma$ is not stable anymore. This scenario is illustrated in Fig.~\ref{fig:time_freezing_2d_exit}. We conclude that under these conditions $y$ leaves $\Sigma$ and enters $R_1$ (into the green region with $f_c(q)>0$).
In this case $\theta = (1,0)$, hence $y' = f_1(y,u) = (f_{\mathrm{ODE}}(x,u),1)$, which matches the unconstrained CLS dynamics \eqref{eq:cls_free_flight} augmented by the clock state.

\begin{figure}[b]
	{\includegraphics[scale=0.70]{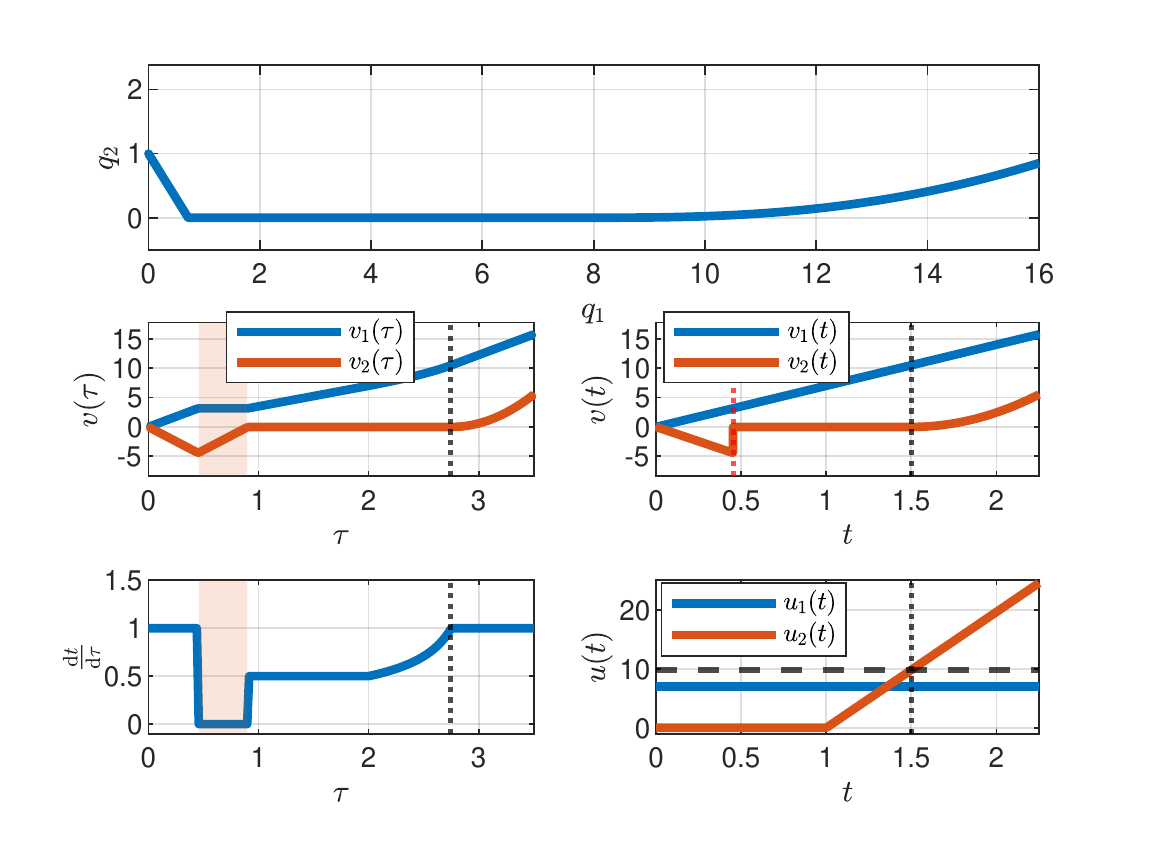}}
	\vspace{-0.7cm}
	\caption{Trajectories of the time-freezing system corresponding to the example CLS. The top plot shows the position of the particle. The middle left plot shows the continuous velocities $v_1$ and $v_2$ in numerical time $\tau$ and the middle right plot shows the discontinuous velocities $v_1$ and $v_2$ in physical time $t$. The bottom left plot shows the speed of time $\frac{\dd t}{\dd \tau}$ and the bottom right plot the control function $u(t)$.
 }
	\label{fig:time_freezing_guiding_example1}
\end{figure}

Effectively, the time-freezing system switches between the DAE and ODE modes, just as the CLS and the state jump is performed by the auxiliary dynamics while the time is frozen. For $\varphi(x,u) \leq 0$ (persistent contact) the solution of the time-freezing system stays on $\Sigma$, just as the solution of the CLS. It leaves the sliding mode when $\varphi(x,u) >0$, which corresponds to contact breaking in the CLS. This relationship is formalized in the next subsection. To illustrate the developments of the last two subsections, we revisit our guiding example and provide a simulation that encompasses all effects discussed so far.
\begin{example}(Speed of time and sliding modes)\label{ex:guiding_example_simulation}
Let us consider the time-freezing PSS from Example \ref{ex:guiding_example_pss}. We choose $a_n = g$. 
It follows that $D(q) = 1$ and $\varphi(x,u) = -g+u_2$. 
We choose a control function 
\begin{align*}
	u(t) = \begin{cases}
		(7,0), & t<1,\\
		(7,2g(t(\tau)-1)), & t\geq 1.  
	\end{cases}
\end{align*}
\noindent Let us make a simulation of the time-freezing system with $y(0) = (0,1,0,0,0)$ for $\tau \in [0,3.5]$. 
The result is depicted in Fig.~\ref{fig:time_freezing_guiding_example1}. 
The particle hits the ground, slides horizontally on it, and lifts off when the control force $u_2(t)$ is stronger than gravity, cf. the top plot. 
We see that when the particle hits the ground, the time is frozen and the auxiliary ODE is active (red strips). 
The vertical velocity $v_2$ becomes zero with the rate $\anormal$.
The system is then a sliding mode with a time slow down factor of $\gamma(x,u) = \frac{\anormal}{\anormal-(-g+u_2)}=\frac{g}{g-(-g+0)}=0.5$ for $t<1$ ($\tau <2$), cf. bottom left plot.
At $t=1$, which corresponds to $\tau=2$ the vertical control force becomes nonzero and $\gamma(x,u)=\frac{g}{2g+u_2}$ grows. For $\tau > 2.8$ we have $\varphi(x(\tau),u(\tau))= -g+u_2 > 0$ and $\Sigma$ is not stable anymore.
The particle lifts off and the contact breaks. 
Note that the solution of the time-freezing system $y(\cdot)$ is continuous in numerical time $\tau$ (middle left plot) and discontinuous in physical time $t$ (middle right plot).
\end{example}

\subsection{Solution relationship}\label{sec:solution_relationship}
We formalize now how to recover the solution of the initial value problem corresponding to the CLS \eqref{eq:cls} from the solution of the time-freezing system from Definition~\ref{def:time_freezing_system}.
\begin{theorem}[Solution relationship]\label{th:solution_relationship}
Regard the initial value problems corresponding to: 
	i) the time-freezing system in Eq. \eqref{eq:time_freezing_system} with a given $y(0) = (q_0,v_0,0)\in \R^{n_y}$ and ${f_c(q_0) \geq 0}$ on a time interval $[0,\tau_{\mathrm{f}}]$, and 
	ii) the CLS from Eq. \eqref{eq:cls} with the initial value $x(0)=(q_0,v_0)\in \R^{n_x}$ on a time interval $[0,\tf]\coloneqq[0,t(\tauf)]$, with $f_c(q(\tf)) \geq 0$ and $n(q(\tf))^\top{v}(\tf)\geq0$. 
	Suppose the following assumptions hold:
	\begin{enumerate}[(a)]
		\item the auxiliary dynamics $f_{\mathrm{aux,n}}(y)$ from Proposition \ref{prop:aux_forming_ode} is used in the time-freezing system in Definition \ref{def:time_freezing_system},
		\item there is at most one time point $\ts=t(\taus)$ where $f_c(q(\ts)) = 0$ and $n(q(\ts))^\top v(\ts^-)<0$ on the time interval $[0,\tf]$,
	\end{enumerate}
	Then, the solutions to the two problems are related as follows:
	\begin{enumerate}
		\item for $t \neq \ts$:
		\begin{subequations}
			\begin{align}
				x(t(\tau)) &= R y(\tau), 
				\text{ with }
				R = \begin{bmatrix}
					{I}_{n_x}  &\mathbf{0}_{n_x,1}
				\end{bmatrix},
				\label{eq:SolutionRelation1_diff}
				\\
				\zn(t(\tau)) & 				= \begin{cases}
				\!-\! D(q(t(\tau)))\varphi(x(t(\tau)), & \mathrm{if}\; y(\tau)\in \Sigma,\\
				0,	 & \mathrm{otherwise}.
				\end{cases} \label{eq:SolutionRelation1_alg} 
			\end{align}
		\end{subequations}
		\item for $t = \ts$:
		\begin{align}\label{eq:SolutionRelation2}
		 \lim_{\substack{\epsilon \to 0\\ \epsilon>0}}	\int_{t_{\mathrm{s}}-\epsilon}^{t_{\mathrm{s}}+\epsilon} \zn(t) \dd t = \int_{\taus}^{\taur} \anormal \dd \tau.
		\end{align}
	\end{enumerate}
\end{theorem}
\textbf{PROOF.} See Appendix \ref{app:proof_of_th1}. \qed

Equivalence for several subsequent impacts is trivially obtained by sequentially applying the argument of the last theorem.
Time-freezing enables one to make the state jump in "slow motion". By plotting the state as a function of physical time we make the "slow" transition "infinitely fast" and recover the discontinuity in time.
More formally, this is encapsulated in \eqref{eq:SolutionRelation2}, which shows that the integral of a Dirac impulse $\zn(t)$ is the same as the integral of the $v$-state of the auxiliary ODE over a finite time interval $[\taus,\taur]$ of nonzero length. 
To simulate a time-freezing system with Zeno's effect, the numerical time horizon would have to be infinitely long, as every state jump requires $\tau_{\mathrm{jump}}>0$. 
In this paper, we assume to have a finite number of impacts.
In practical robotics applications, one is usually interested in solutions with a finite number of impacts.

\subsection{Possible extensions}
In this paper, we consider a single unilateral constraint. To extend the same ideas for multiple and/or simultaneous impacts one must take care of the chosen impact model \cite{Brogliato2016}.
The extension can be made in several ways, e.g., time continues to flow when the first normal velocity component reaches zero or when all of them reach zero. Some multiple impact models suffer from nonuniqueness of solutions \cite{Brogliato2016,Nguyen2018,Stewart2000}, and how to proceed is a modeling decision.
In the special case when the constraints, e.g., $f_{c,i}(q)$ and $f_{c,j}(q)$ are orthogonal in the kinetic metric, i.e., $f_{c,i}(q)^\top M(q)^{-1} f_{c,j}(q) = 0$, the impacts can be treated independently \cite{Brogliato2016}.
In this case, one would take for every constraint the auxiliary dynamics from Prop. \ref{prop:aux_forming_ode} and the time would continue to flow when all normal velocity components reach zero. 
Stochastic approaches with multiple outcomes are also possible, e.g., within the binary collision law as shown in \cite{Nguyen2018} or a stochastic version of Routh's impact model~\cite{Halm2021}.

\section{Frictional impact}\label{sec:frictional_impact}
If friction is present at the contact point, frictional impulses cause state jumps in the tangential directions.
This section extends the time-freezing reformulation for CLS with frictional impacts.  Appropriate auxiliary dynamics for the state jumps in the tangential directions are introduced. The time-freezing system covers both stick and slip motions.
\subsection{CLS with frictional impacts}\label{sec:cls_friction}
The extension of the CLS \eqref{eq:cls} with Coulomb friction reads as:
\begin{subequations}\label{eq:cls_friction}
\begin{align}
	\dot{q} &= v, \\
	\dot{v} &= f_\mathrm{v}(q,v,u)  + M(q)^{-1} (n(q) \zn+\Tan\zt),\\
	0 &\leq \zn \perp f_c(q) \geq 0, \\
	\begin{split}
		0&=n(q(\ts))^\top v(\ts^+),\\
		&\mathrm{if} \ f_c(q(\ts)) = 0\ \mathrm{and}\ n(q(\ts))^\top v(\ts^-)<0,
	\end{split}\\
	&\zt  \in\arg\min_{{\zttilde} \in \R^{\nt}} \quad  -v^\top \Tan{\zttilde}\label{eq:cls_friction_disipiation1}
	\\ 
	&\quad	\textrm{s.t.} \quad  \| {\zttilde} \|_2 \leq \mu \zn. \label{eq:cls_friction_disipiation2}
\end{align}
\end{subequations}
Compared to the CLS \eqref{eq:cls}, the model is extended by the term $M(q)^{-1}\Tan\zt$ in the r.h.s. of velocity dynamics (friction force) and Eq. \eqref{eq:cls_friction_disipiation1}-\eqref{eq:cls_friction_disipiation2} (friction model).
The matrix $\Tan = [b_1(q),\ldots, b_{\nt}(q)] \in \R^{n_q \times \nt}$ spans the tangent space at contact points $\{q\in \R^{n_q} \mid f_c(q) =0\}$ and $\zt\in \R^{\nt}$ is the friction force. For planar contacts, we have $\nt = 1$, and for 3D contacts $\nt =2$. Denote the tangential velocity during contact by $\vt \coloneqq \Tan^\top v \in \R^{\nt}$.
We assume that the vectors $b_1(q),b_2(q)$ and $n(q)$ are orthogonal in the kinetic metric,
i.e., $b_1(q)^\top M(q)^{-1} n(q) = 0$, $b_2(q)^\top M(q)^{-1} n(q) = 0$ and $b_1(q)^\top M(q)^{-1} b_2(q) = 0$.
While restrictive, this assumption simplifies computations and maintains brevity in the exposition. 
Following a similar approach without it leads to the same conclusions but with more complex equations.
The convex optimization problem \eqref{eq:cls_friction_disipiation1}-\eqref{eq:cls_friction_disipiation2} is the \textit{maximum dissipation principle} \cite{Moreau1977a} and $\mu >0$ is the coefficient of friction. This model expresses that the dissipation of the kinetic energy between two objects in contact is maximized. 
When impacts occur, the impulsive $\zn$ results in an impulsive $\zt$ via Eq. \eqref{eq:cls_friction_disipiation1}-\eqref{eq:cls_friction_disipiation2} and thus we have state jumps in the velocity $v$ in the tangential directions $b_j(q)$, $j =1,\ldots, \nt$. 

\subsection{Stick and slip dynamics of the CLS}\label{sec:cls_slip_stick}

For a given $\zn$ the solution map of the convex optimization problem ~\eqref{eq:cls_friction_disipiation1}-\eqref{eq:cls_friction_disipiation2} is given by \cite{Moreau1977a}:
\begin{align}\label{eq:friction_slution_map}
	\zt &\in  \begin{cases}
		\{- \mu \zn  \frac{ \vt }{\| \vt  \|_2} \}, & \textrm{if}\; \| v_{\mathrm{t}} \|_2 > 0,\\
		\{ {\zttilde} \mid \|\zttilde\|_2 \leq \mu \zn \}, & \textrm{if}\; \| v_{\mathrm{t}} \|_2 = 0.
	\end{cases} 
\end{align}
When the system is in contact, it can be in \textit{slipping motion}, i.e., it has nonzero tangential velocity $\vt \neq 0$, or in \textit{sticking motion} with $ \vt =0$. 
Similar to Section \ref{sec:cls_modes}, we derive equivalent ODEs which model the stick and slip dynamics during contact phases. 
If the system is in slip motion and $\zn>0$, it follows from \eqref{eq:friction_slution_map} that $\zt = - \mu \zn \frac{\vt}{\|\vt \|_2}$. Similar to Eq. \eqref{eq:cls_index_reduced_dae} we obtain
\begin{subequations}\label{eq:cls_slip}
	\begin{align}
	\dot{q} &= v,\\
	\begin{split}
	\dot{v} & = f_{\mathrm{v}}(q,v,u)-M(q)^{-1} D(q)^{-1}\varphi(x,u)  \\
	&\Big(n(q)- \Tan \mu\frac{\vt}{\|\vt \|_2}\Big).
	\end{split}
	\end{align} 
\end{subequations}
The r.h.s. of this ODE is compactly denoted by $\fslip(x,u)$.
If the system is in sticking motion we have $\vt = \Tan^\top v = 0$ and $n^\top v = 0$. By differentiating these equations  w.r.t. time we can explicitly compute the multipliers $\lambda \coloneqq  (\zn,\zt)$ by
\begin{align*}
\lambda &= -\tilde{D}(q)^{-1}   \tilde{\varphi}(x,u), 
\\
\tilde{D}(q)\! &= \! \begin{bmatrix}	n(q)& \! \Tan  \end{bmatrix}^\top 
				\! M(q)^{-\! 1} 
	\begin{bmatrix} n(q)&\! \Tan \end{bmatrix},\\
\tilde{\varphi}(x,u)\! &\coloneqq \!	\begin{bmatrix} n(q)& \! \Tan \end{bmatrix}^\top \! \!	 f_v(q,v,u)\!	 +\!	 
\nabla_q( \begin{bmatrix}	n(q)& \!  \Tan  \end{bmatrix}^\top\!	 v )\!	^\top \!  v. 
\end{align*}	
The ODE describing the sticking motion reads as
\begin{subequations}\label{eq:cls_stick}
	\begin{align}
		\dot{q} &= v,\\
		\dot{v} & = f_{\mathrm{v}}(q,v,u) 
			-M(q)^{-1} \begin{bmatrix}	n(q)& \Tan  \end{bmatrix} \tilde D(q)^{-1}\tilde{\varphi}(x,u).
	\end{align} 
\end{subequations}
Its r.h.s. is compactly denoted by~$\fstick(x,u)$.
The transition from the stick to the slip mode occurs when the other tangential forces are greater than the maximal friction force $\zt$, cf.~\cite[Chapter 5]{Brogliato2016}.
\subsection{Time-freezing for frictional impacts in the 2D case}\label{sec:time_freezing_friction}
In the planar case we have $\nt = 1$ and Eq. \eqref{eq:friction_slution_map} simplifies to $\zt \in -\mu \zn \sign(v_t)$. Denote the single column of $\Tan$ by $b(q)$. Depending on the sign of $\vt$, we define an auxiliary dynamical system to mimic the state jump in the tangential direction $b(q)$. For the $n(q)$-direction we use the dynamics from Proposition \ref{prop:aux_forming_ode}. For the tangential direction and $b(q)^\top v < 0$ we define the \textit{tangential auxiliary dynamics} analogously:
\begin{align} \label{eq:auxiliary_dynamics_tangetnial}
	y'=f_{\mathrm{aux,t}}^-({y}) &
	\coloneqq \begin{bmatrix}
		\mathbf{0}_{n_q , 1} \\   M({q})^{-1} 	b(q) \atangent \\0
	\end{bmatrix}.
\end{align} 
To account for the sign of the tangential velocity, for $b(q)^\top v > 0$ we use $y'=f_{\mathrm{aux,t}}^+(y) \coloneqq -f_{\mathrm{aux,t}}^-(y)$. Depending on the sign of $\vt$, one of these ODE is active for the same numerical time interval of the length $\taujump$ as $y' = f_{\mathrm{aux,n}}(y)$. Furthermore, we know from  Eq.~\eqref{eq:SolutionRelation2} in Theorem~\ref{th:solution_relationship} that the impulse bringing the normal velocity $n(q)^\top v<0$ to zero after an impact is proportional to $\anormal\taujump$. Thus, by settings $\atangent = \mu \anormal$, we conclude that the integrals of the auxiliary dynamics satisfy the maximum dissipation principle, i.e., $\atangent  \taujump =  - \mu \anormal\taujump  \mathrm{sign} (b(q)^\top v)$. 

State jumps in both the normal and tangential directions are treated simultaneously with different auxiliary dynamics. 
They should be active whenever $y \in Q \coloneqq \{y \in \R^{n_y} \mid  c_1(y) <0 ,c_2(y)<0 \}$, cf. Fig. \ref{fig:time_freezing_2d}. To treat different signs of the tangential velocity we introduce the switching function $ c_3(y) = b(q)^\top v$.
Hence, we have in total $n_f = 3$ regions, one for the unconstrained dynamics and two to mimic the state jumps. 
We extend the definition of the regions in Eq. \eqref{eq:time_freezing_regions} as follows:
\begin{subequations}\label{eq:time_freezing_regions_friction}
\begin{align}
		\begin{split}
		R_{1} &	=  \{y \in \R^{n_y} \mid c_{1}(y)>0\}\\
		& \cup \{y \in \R^{n_y} \mid c_{1}(y)<0, c_{2}(y)>0 \},\ 
		\end{split}\\
		R_{2} & = Q \cap \{y \in \R^{n_y} \mid c_{3}(y)>0\},\\
		R_{3} & = Q \cap \{y \in \R^{n_y} \mid c_{3}(y)<0\}.
	\end{align}
\end{subequations}
The sum of the corresponding auxiliary dynamics accounts for the simultaneous state jumps, i.e.,
$f_2(y)  =\! f_{\mathrm{aux,n}}(y)\! + \! f_{\mathrm{aux,t}}^+(y),\; f_3(y) =\! f_{\mathrm{aux,n}}(y)\! +\! f_{\mathrm{aux,t}}^-(y)$.
The time-freezing system for the CLS \eqref{eq:cls_friction} is given in the next definition.
\begin{definition}[Time-freezing system with friction]\label{def:time_feezing_friction}
Let $\tau \in \R$ be the numerical time, $y(\tau) \coloneqq ({x}(\tau),t(\tau)) \in \R^{n_y}$ the differential states and $u(\tau)\in\R^{n_u}$ a given control function. 
The time-freezing PSS is defined by the regions $R_{i}, i =1,\ldots,3$, in  Eq. \eqref{eq:time_freezing_regions_friction}
with 
\begin{align*}
f_{1}(y,u) &= (f_{\mathrm{ODE}}(x,u),1),\\
f_{2}(y) &= f_{\mathrm{aux,n}}(y)+f_{\mathrm{aux,t}}^+(y),\\
f_{3}(y) &= f_{\mathrm{aux,n}}(y)+f_{\mathrm{aux,t}}^-(y).
\end{align*}
The corresponding Filippov system, which we call the time-freezing system, is denoted by $y'  \in F_{\mathrm{TF}}(y,u)$. 
The set $F_{\mathrm{TF}}(y,u)$ is defined as in \eqref{eq:FilippovDI_with_multiplers}.
It is assumed that appropriate dynamics $f_{\mathrm{aux,n}}({y}),f_{\mathrm{aux,t}}^+(y)$ and $f_{\mathrm{aux,t}}^-(y)$ exist.
\end{definition}
As in the frictionless case, we are interested in the relation of the CLS in contact mode and the corresponding sliding mode of the time-freezing system on $\Sigma$.
\begin{theorem}[Slip-stick sliding mode]\label{th:slip_stick_theorem}
Suppose that the auxiliary dynamics from Proposition \ref{prop:aux_forming_ode} and Eq. \eqref{eq:auxiliary_dynamics_tangetnial} are used in the time-freezing system from Definition~\ref{def:time_feezing_friction}. Let $y(\tau)$ be a solution of this system with $y(0) \in \Sigma$ and $\tau \in [0,\tauf]$. 
Suppose that $\varphi(x(\tau),u(\tau)) \leq 0$ for all $\tau \in [0,\tauf]$ (persistent contact), then the following statements are true:
\vspace{-0.35cm}
\begin{enumerate}[(i)]
	\item If $\vt \neq 0$ (slip motion), then the sliding mode dynamics are given by ${y'=\gamma(x,u)(\fslip(x,u),1)}$.
	\item If $\vt =0$ (stick motion), then the sliding mode  dynamics are given by ${y'=\gamma(x,u)(\fstick(x,u),1)}$,  
\end{enumerate}
where $\gamma(x,u)\in \left(0,1\right]$ is a time-rescaling factor defined in Eq. \eqref{eq:time_freezing_rescaling}.
\end{theorem}
\textbf{PROOF.} See Appendix \ref{app:proof_of_th2}. \qed

This result generalizes Theorem \ref{th:presistent_contact} and one can see that the sliding mode dynamics match the slip or stick dynamics of the CLS. Finally, we show that CLS with planar contacts, friction, and impacts are equivalent to Filippov systems.
\begin{theorem}[Solution relationship]\label{th:solution_relationship_friction}
	Regard the initial value problems corresponding to: 
	i) the time-freezing system in Definition \ref{def:time_feezing_friction} with a given $y(0) = (q_0,v_0,0)\in \R^{n_y}$ and ${f_c(q_0) \geq 0}$ on a time interval $[0,\tau_{\mathrm{f}}]$, 
	ii) the CLS from Eq. \eqref{eq:cls_friction} with the initial value $x(0)=(q_0,v_0)\in \R^{n_x}$ on a time interval $[0,\tf]\coloneqq[0,t(\tauf)]$, with $f_c(q(\tf)) \geq 0$ and $n(q(\tf))^\top{v}(\tf)\geq0$. 
	Suppose the following assumptions hold:
	\begin{enumerate}[(a)]
		\item the auxiliary dynamics $f_{\mathrm{aux,n}}(y)$ from Proposition \ref{prop:aux_forming_ode} and 
		$f_{\mathrm{aux},t}^-(y)$, $f_{\mathrm{aux},t}^+(y)$ from Eq. \eqref{eq:auxiliary_dynamics_tangetnial}
		 are used in the time-freezing system in Definition \ref{def:time_feezing_friction},
		\item there is at most one time point $\ts=t(\taus)$ where $f_c(q(\ts)) = 0$ and $n(q(\ts))^\top v(\ts^-)<0$ on the time interval $[0,\tf]$,
	\end{enumerate}
	Then, the solutions to the two problems are related as follows:
	\begin{enumerate}
		\item For $t \neq \ts$:
			\begin{align*}
				x(t(\tau)) &= R y(\tau), 
				\text{ with }
				R = \begin{bmatrix}
					{I}_{n_x}  &\mathbf{0}_{n_x,1}
				\end{bmatrix},
			\\
				\lambda(t(\tau)) & = \begin{cases}
					\lambda_{\mathrm{Slip}}(t(\tau)), & \mathrm{if}\; \ y\in \Sigma,\; \vt \neq 0,\\
					\lambda_{\mathrm{Stick}}(t(\tau)), & \mathrm{if}\; \ y\in \Sigma,\; \vt  = 0,\\
					0,	 & \mathrm{otherwise}.
				\end{cases} 
		\end{align*}
		\vspace{-0.25cm}
		with 
		\vspace{-0.15cm}
		\begin{align*}
		\lambda_{\mathrm{Slip}} &=-{D}(q)^{-1}\varphi(x,u) (1,-\mu \sign(\vt))\\
		\lambda_{\mathrm{Stick}} &= -\tilde{D}(q)^{-1}\tilde{\varphi}(x,u).
		\end{align*}
		\item For $t = \ts$:
		\label{eq:solution_relationship_weak}
		\begin{align*}
			\lim_{\substack{\epsilon \to 0\\ \epsilon>0}}	\int_{t_{\mathrm{s}}-\epsilon}^{t_{\mathrm{s}}+\epsilon} \zn(t) \dd t &= \int_{\taus}^{\taur} \anormal \dd \tau,\\
			\lim_{\substack{\epsilon \to 0\\ \epsilon>0}}	
			\int_{t_{\mathrm{s}}-\epsilon}^{t_{\mathrm{s}}+\epsilon} |\zt(t)| \dd t &=  \int_{\taus}^{\taur} 
			\mu\anormal \dd \tau,
		\end{align*}
	\end{enumerate}
\end{theorem}
\textbf{PROOF.}  Theorem \ref{th:slip_stick_theorem} is applied for $y\in \Sigma$, otherwise the proof follows similar lines as the proof of Theorem \ref{th:solution_relationship}. The absolute value in the last equation accounts for all signs of $\vt$. \qed 

\begin{figure}[t]
	\centering
	{\includegraphics[scale=0.65]{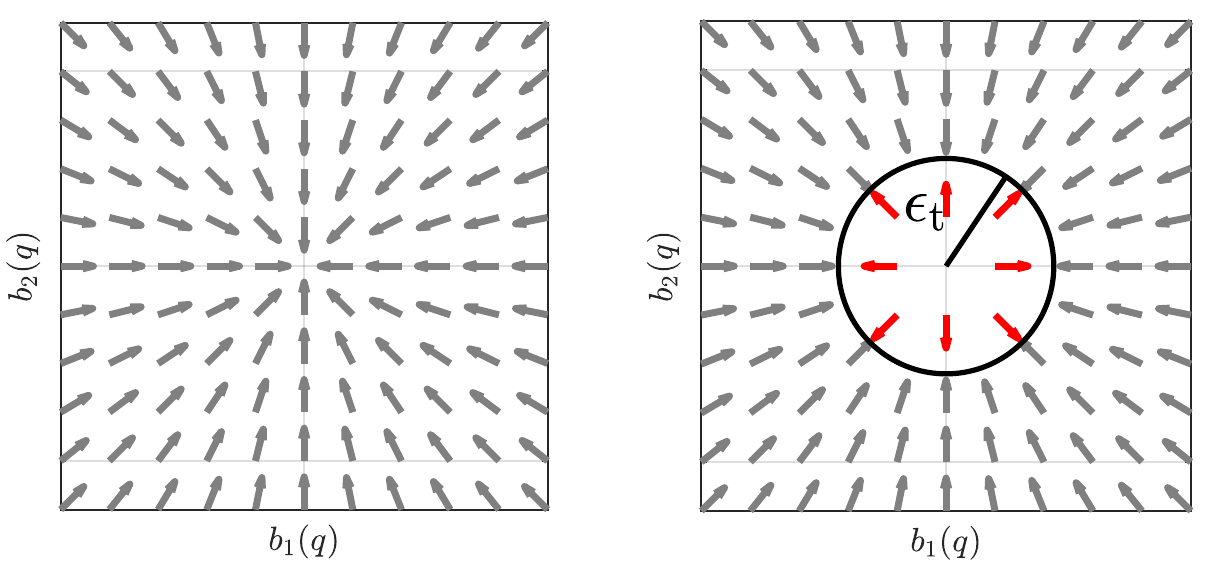}}
	\vspace{-0.5cm}
	\caption{The exact friction force Eq. \eqref{eq:friction_slution_map} (left) and its approximation in Eq. \eqref{eq:friction_slution_map_approx} (right). }
	\label{fig:friction_relaxed}
\end{figure}

\subsection{Time-freezing for frictional impacts in the 3D case}
This case is more difficult since we cannot easily treat different directions of $\vt$ with different auxiliary dynamics as in the planar case. The solution map of \eqref{eq:friction_slution_map} depends discontinuously on $\|\vt\|_2$, hence, we must take it as a switching function. The set $\| \vt\|_2=0$ has no interior and we cannot use the Filippov extension from Eq.~\eqref{eq:FilippovDI_with_multiplers}, which assumed regions $R_i$ with nonempty interior. More general definitions without multipliers $\theta$ that could treat this case exist \cite{Filippov1988}, but they are not computationally useful in our case, as we see in the next section. To alleviate this difficulty we propose an approximation for~\eqref{eq:friction_slution_map}:
\begin{align}\label{eq:friction_slution_map_approx}
	\zt & =  \begin{cases}
		- \mu \zn  \frac{ v_{\mathrm{t}} }{\| v_{\mathrm{t}}  \|_2}, & \textrm{if}\; \| v_{\mathrm{t}} \|_2 > \epst,\\
		{ v_{\mathrm{t}}}, & \textrm{if}\; \| v_{\mathrm{t}} \|_2 < \epst,
	\end{cases} 
\end{align}
with a small parameter $\epst>0$. This expression is exact for $\| \vt \|_2 > \epsilon_{\mathrm{t}}$, thus we can make it arbitrarily accurate. For $\| \vt \|_2 < \epsilon_{\mathrm{t}}$, the vector field drives the tangential velocity towards $\| \vt \|_2 = \epsilon_{\mathrm{t}}$, see Fig.~\ref{fig:friction_relaxed}. In a Filippov setting, a convex combination of the two cases in \eqref{eq:friction_slution_map_approx} keeps the velocity at $\| \vt \|_2 = \epsilon_{\mathrm{t}}$. 
Hence, in sticking mode, we have a velocity drift of $\epst$.
Now by taking $c_3(y) = \| \vt \|_2 - \epst$ as a switching function, we can define regions with nonempty interiors and the corresponding auxiliary dynamics.
The auxiliary dynamics mimicking the behavior of \eqref{eq:friction_slution_map_approx} read as
\begin{align*}
		f_{\mathrm{aux,t}}^+(y)\! &= \! \begin{bmatrix}
			\mathbf{0}_{n_q , 1} \\   -M(q)^{-\!1} \!	B(q) \atangent  \frac{\vt}{\| \vt\|}\\0
		\end{bmatrix},\\
		f_{\mathrm{aux,t}}^-(y)\!&= \! \begin{bmatrix}
		\mathbf{0}_{n_q , 1} \\  M(q)^{-1}B(q) {\vt}\\0
	\end{bmatrix}.
\end{align*} 
\noindent The regions for the time-freezing system are defined as in \eqref{eq:time_freezing_regions_friction} and the matching time-freezing system is defined analogously to Definition~\ref{def:time_feezing_friction}.
Furthermore, one could derive stick-slip dynamics corresponding to the solution map approximation \eqref{eq:friction_slution_map_approx} and relate it to the time-freezing system by following similar lines as in Theorems \ref{th:slip_stick_theorem} and \ref{th:solution_relationship_friction}, but we omit the details here.
We conclude this section by revisiting Example \ref{ex:guiding_example_simulation}, but now with adding friction.
\begin{figure}[t]
	{\includegraphics[scale=0.70]{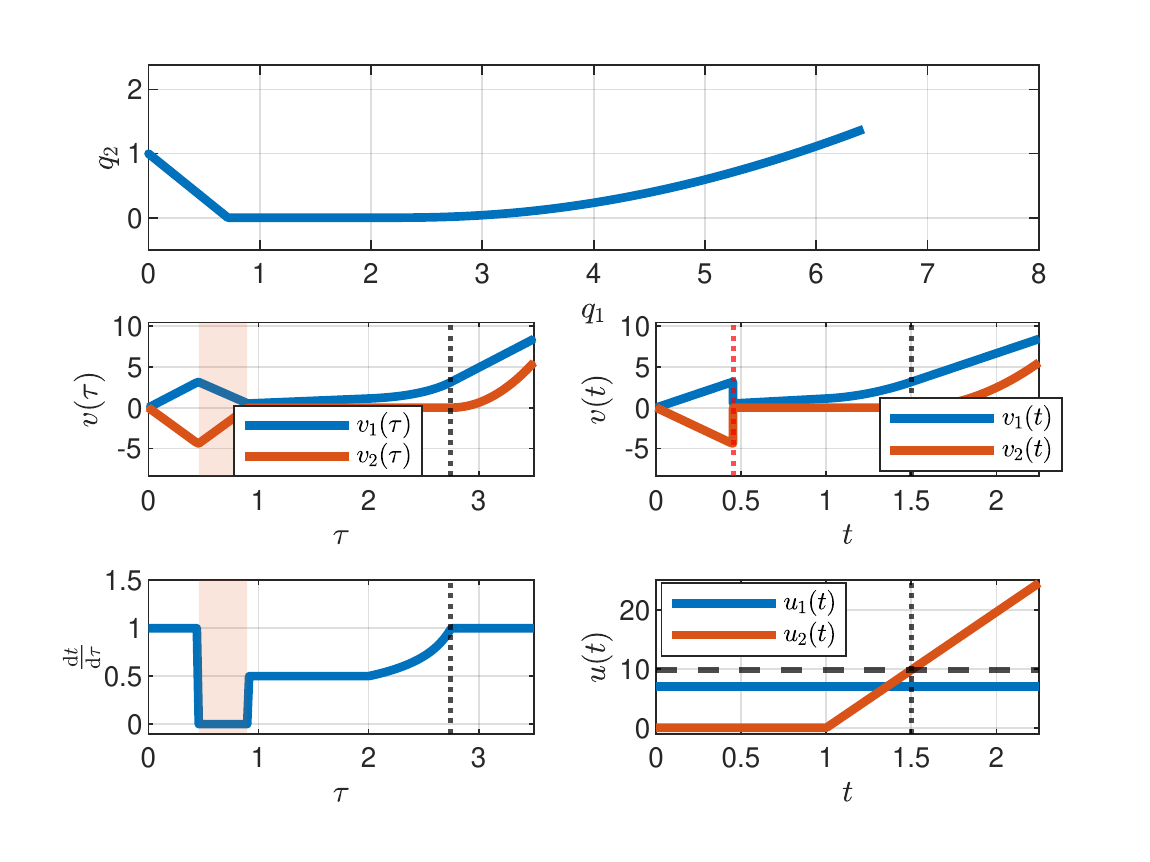}}
	\vspace{-0.6cm}
	\caption{Trajectories of the time-freezing system from Ex.~\ref{ex:guiding_example_simulation_friction}.}
	\label{fig:time_freezing_motivation2}
\end{figure}
\begin{example}(Frictional impact)\label{ex:guiding_example_simulation_friction}
The time-freezing system model from Example \ref{ex:guiding_example_pss} is extended by adding friction with a coefficient $\mu = 0.6$. 
In the planar case, we have $\nt =1$ and the tangent at the contact point is 
$b(q) = [1,0]$.  
We have the switching functions $c(y) = (q_2,v_2,v_1)$.
Following Eq.~\eqref{eq:time_freezing_regions_friction}, the regions of the time-freezing system are 
$R_1 = \{ y\mid q_2>0 \} \cup \{y \mid q_2 <0, v_2>0\}$,
$R_2 =\{y \mid q_2 <0, v_2<0 , v_1>0\}$ and
$R_3 =\{y \mid q_2 <0, v_2<0 , v_1<0\}$. 
The dynamics of the PSS are
$f_1 = (v_1,v_2,u_1,-g+u_2,1)$,
$f_2 = (0,0,-\mu\anormal,\anormal,0)$ and
$f_3 = (0,0,\mu\anormal,\anormal,0)$. 
The results of the simulation are depicted in Fig.~\ref{fig:time_freezing_motivation2}. 
Note that due to friction there is now also a state jump in the tangential velocity $v_1$, cf. middle plots. 
Afterward, the acceleration of $v_1$ is during contact phases smaller due to the friction force.
However, the tangential acceleration is increasing over time as the normal contact force becomes weaker because of $u_2$. 
At $\tau=2.8$ the particle lifts off as in the previous example. 
\end{example}

\section{The application of time-freezing in optimal control problems}\label{sec:computational_considerations}
This section regards OCP formulations with time-freezing systems and numerical methods to solve them. 
We derive an equivalent OCP, now subject to the time-freezing system. 
The section concludes with discussing the numerical methods and software for solving OCPs subject to time-freezing systems.
\subsection{Continuous-time OCP with a CLS}
We regard a modification of the OCP \eqref{eq:ocp_cls1}, where we consider the CLS with friction \eqref{eq:cls_friction}.
This continuous-time OCP read as:

\begin{subequations} \label{eq:ocp_cls}
	\begin{align}
		\min_{x(\cdot),\lambda(\cdot),u(\cdot),} \quad &   \Psi(x(\Tctrl)) \label{eq:ocp_cls_objective}\\
		\textrm{s.t.} \quad  &x(0) = \bar{x}_0 \label{eq:ocp_cls_iv},\\
		&\textrm{Eq.} \eqref{eq:cls_friction},\; t \in [0,\Tctrl] \label{eq:ocp_cls_dynamics}\\
		&0\leq g(x(t),u(t)),\; t \in [0,\Tctrl] \label{eq:ocp_cls_path},\\
		&0\leq r(x(\Tctrl)),\;  \label{eq:ocp_cls_temrinal}
	\end{align}
\end{subequations}
The functions $g: \R^{n_x} \times \R^{n_u} \to \R^{n_g}$ and $r: \R^{n_x} \to \R^{n_r}$ are the path and terminal constraints, respectively. 
W.l.o.g. we only consider a terminal cost term here and remind the reader that the integral of a running cost $L: \R^{n_x} \times \R^{n_u} \to \R$ over $[0,\Tctrl]$ can be treated via a terminal cost term by introducing a quadrature state~$\ell(t)$
\begin{align}\label{eq:quadrature_cls}
	\frac{\dd }{\dd t}  \ell(t) & = L(x(t),u(t)), t\in [0,\Tctrl],\; 	\ell(0)  =0,
\end{align}
and adding $\ell(T)$ to the objective.

\subsection{Continuous-time OCP with a time-freezing system}
Based on the results from the previous sections, we derive now an OCP subject to a time-freezing system.
The new OCP is regarded in numerical time $\tau \in [0,\tilde{\Tctrl}]$.
We take four steps in this transformation:
(1) we modify the quadrature state in Eq. \eqref{eq:quadrature_cls} so that the cost integrated over both numerical and physical time remains unchanged;
(2) we reformulate the time-freezing system into an equivalent dynamic complementarity system to make it possible to apply FESD;
(3) we introduce a time-transformation to ensure that the terminal physical time $t(\tilde{\Tctrl})$ matches the true control horizon of Eq.\eqref{eq:ocp_cls}, i.e., $t(\tilde{\Tctrl}) = \Tctrl$;
(4) we express the remaining constraints in terms of numerical time.

We start with adapting the objective.
This is achieved by replacing the quadrature state \eqref{eq:quadrature_cls} by:
\begin{align}\label{eq:quadrature_tf}
	\frac{\dd}{\dd \tau }\ell(\tau) &= \begin{cases}
		L(x(\tau),u(\tau)),&\; \textrm{if}\; y \in R_1,\\
		0,&\; \textrm{otherwise}.
	\end{cases}
\end{align}
When the time is frozen the cost integral is zero and there are no contributions to the overall objective, i.e.,  the cost is unchanged when the time is frozen.

Next, we rewrite the time-freezing system from Definition \ref{def:time_feezing_friction} as an equivalent dynamic complementarity system.
This enables the application of the FESD method to the time-freezing system.
To achieve this, we use the set-valued step function $\alpha(x) = 0.5(1+\Sign(x))$. 
It is well-known that Filippov multipliers $\theta$ can be expressed as products of set-valued step functions $\alpha_i$ \cite{Acary2014,Dieci2011}.
We define $c(y) \coloneqq (c_1(y),c_2(y),c_3(y))$.
The set-valued step function $\alpha(\cdot)$ can be expressed as the solution map of a parametric linear program \cite{Acary2014}:
\begin{align}\label{eq:parametric_lp}
	\alpha(c(y)) \in \arg\min_{\tilde{\alpha}\in \R^3} -c(y)\tilde{\alpha} \; \mathrm{s.t.} \; 0\leq \tilde{\alpha} \leq e. 
\end{align}
Using the Karush–Kuhn–Tucker (KKT) conditions of \eqref{eq:parametric_lp} we can derive from $y' \in F_{\mathrm{TF}}(y,u)$ in Definition \ref{def:time_feezing_friction} the equivalent dynamic complementarity system:
\begin{subequations}\label{eq:time_freezing_dcs}
	\begin{align}
		y' &= F(y,u)~\theta,\\
		0 & = g_{\mathrm{F}}(\theta,\alpha),\\
		0 &= c(y) - \lambdap  + \lambdan,\\
		0 &\leq \alpha \perp \lambdan \geq 0,\\
		0 &\leq e-\alpha \perp \lambdap \geq 0.
	\end{align}
\end{subequations}
The matrix $F(y,u) =[f_1(y,u),\ldots f_{n_f}(y)] \in \R^{n_y \times n_f}$ collects the modes of the PSS and $\theta = (\theta_{1},\ldots,\theta_{n_f})$. 
The last three lines are the KKT conditions of \eqref{eq:parametric_lp}, where $\lambdan,\lambdap\in \R^3$ are the Lagrange multipliers for the lower and upper bounds in \eqref{eq:parametric_lp}, respectively. 
We group all algebraic variables of the DCS in the vector $z = (\theta,\alpha,\lambdap,\lambdan)$.
The function $g_{\mathrm{F}}$ relates the Filippov multipliers $\theta$ with the evaluations of the step functions~$\alpha$:
\begin{align}\label{eq:time_freezing_dcs_lifting}
	g_{\mathrm{F}}(\theta,\alpha) \coloneqq \begin{bmatrix}
		\theta_{1}-\alpha_1 + (1-\alpha_1)\alpha_2,\\
		\theta_{2}-(1-\alpha_1)(1-\alpha_2)(1-\alpha_3)\\
		\theta_{3}-(1-\alpha_1)(1-\alpha_2) (\alpha_3)\\
	\end{bmatrix}.
\end{align}
These expressions correspond to the signs of $c_j(y)$ in the definitions of $R_i$, e.g., $c_1(y)>0, c_2(y)<0$ results in $\alpha_1(1-\alpha_2)$ \cite{Acary2014,Dieci2011}. 

The time-freezing system evolves over $\tau \in [0,\tilde{\Tctrl}]$. 
During state jumps the physical time evolution is stopped.
As a consequence, we have that $t(\tilde{\Tctrl}) < \Tctrl$, i.e., the terminal physical time in the time-freezing problem does not match the desired time $\Tctrl$. 
To resolve this, we introduce a time-transformation variable $s(\tau) \in \R$ and impose the terminal constraint on the clock state $t(\tilde{\Tctrl}) = \Tctrl$. 
The function $s(\cdot)$ can be interpreted as virtual control that controls the physical time in the numerical time.
Consequently, we obtain $t' = s$ and with $s >1$, which speeds up the physical time and allows us to catch up and reach the desired time $\Tctrl$, cf. Example \ref{ex:illustrate_ocp}.
Such time transformations are very common in optimal control when one wants to optimize over the terminal time \cite{Rawlings2017}.

It is left to impose the path \eqref{eq:ocp_cls_path} and terminal constraints \eqref{eq:ocp_cls_temrinal} in numerical time for $x(\tau)$ and $u(\tau)$.
Finally, the OCP subject to the time-freezing system reads as:
\vspace{-0.1cm}
\begin{subequations} \label{eq:ocp_time_freezing}
	\begin{align}
		\min_{\substack{y(\cdot),z(\cdot),\\u(\cdot),s(\cdot)}} \quad &   \Psi(x(\tilde{\Tctrl})) \label{eq:ocp_time_freezing_objective}\\
		\textrm{s.t.} \quad  &x(0) = \bar{x}_0,\; t(0) = 0, \label{eq:ocp_time_freezing_iv}\\
		&y'(\tau)\! =\! s(\tau)F(y(\tau),u(\tau)\!)\theta(\tau), \tau\!\in \![0,\!\tilde{\Tctrl}],\label{eq:ocp_time_freezing_dcs_1}\\
		&0  = g_{\mathrm{F}}(\theta(\tau),\alpha(\tau)), \tau\in[0,\tilde{\Tctrl}],\label{eq:ocp_time_freezing_dcs_theta}\\
		&0\! = \! c(y(\tau)) - \lambdap(\tau)  + \lambdan(\tau), \tau\in[0,\tilde{\Tctrl}],\label{eq:ocp_time_freezing_dcs_switching}\\
		&0 \leq \alpha(\tau) \perp \lambdan(\tau)\geq 0, \tau\in[0,\tilde{\Tctrl}],		\label{eq:ocp_time_freezing_dcs_comp1}\\
		&0 \! \leq\! e \!-\! \alpha(\tau)\! \perp \! \lambdap(\tau) \! \geq\! 0, \tau\! \in \! [0,\tilde{\Tctrl}],\!\label{eq:ocp_time_freezing_dcs_comp2}\\
		&0\leq g(x(\tau),u(\tau)),\; \tau\! \in \! [0,\tilde{\Tctrl}], \label{eq:ocp_time_freezing_path}\\
		&0\leq r(x(\tilde{\Tctrl})),\;  \label{eq:ocp_time_freezing_temrinal}\\
		&t(\tilde{\Tctrl}) = {\Tctrl}.\;  \label{eq:ocp_time_freezing_temrinal_time}
	\end{align}
\end{subequations}
It is important to note, that when the time is frozen ($t' =0$) the control $u(\tau)$ does not influence $x(\tau)$, since the auxiliary dynamics does not depend on the control, cf. Eq \eqref{eq:aux_dyn_example}. 
Therefore, one could even omit the path constraints whenever $t' = 0$, but we keep it for notational simplicity. 
Additionally, the integral of the stage cost remains unchanged, since $\frac{\dd}{\dd \tau}\ell(\tau)= 0$ in this case, cf.~ \eqref{eq:quadrature_tf}.

Next, we show that the optimal controls obtained by solving the initial OCP \eqref{eq:ocp_cls}, with appropriate modifications, are also optimal for \eqref{eq:ocp_time_freezing}.
Let $u^*(t), t\in [0,\Tctrl]$ be an optimal control of \eqref{eq:ocp_cls}. 
We construct an $\tilde{u}^*(\tau)$, $t\in [0,\tilde{\Tctrl}]$ as follows. 
It can be seen that when the physical time is evolving ($t'>0$), we can find the inverse function $t^{-1}(\cdot)$ to find the corresponding numerical time $\tau$. 
We construct a control function for the time-freezing system:
\begin{align}\label{eq:optimal_control}
	\tilde{u}^*(\tau) &= 
	\begin{cases}
		u(t^{-1}(t(\tau))), & \textrm{ for } t(\tau)' > 0 \\
		\hat{u}(\tau), & \textrm{ for } t(\tau)' = 0,
	\end{cases} 
\end{align} 
where $\hat{u}(\tau)$ is any function such that $g(x(\tau),\hat{u}(\tau))\geq0$ holds, whenever $t'(\tau) = 0$. 
Recall that $\hat{u}(\tau)$ does not change the objective nor it changes ${x}(\tau)$, its only purpose is to extend $u(t)$ to intervals when the time is frozen.
For example, we can choose a constant value that does not violate the path constraints.
With \eqref{eq:optimal_control} we can extend $u^*(t)$ for the time interval where the physical time is frozen. 
Conversely, given an optimal control $\tilde{u}^{*}(\tau)$ of Eq. \eqref{eq:ocp_time_freezing}, then
we expect ${u}^{*}(t(\tau))$ to be optimal for~\eqref{eq:ocp_cls}.

\begin{theorem}\label{th:optimal_controls}
	Let $\tilde{u}^*(\tau), \; \tau\in [0,\tilde{\Tctrl}]$ be an optimal control obtained by solving the OCP \eqref{eq:ocp_time_freezing}.
	Then $u^*(t) = \tilde{u}^*(t(\tau)),\; t\in [0,T]$ is an optimal control of the OCP \eqref{eq:ocp_cls}. 
	Conversely, let $u^*(t),\; t\in [0,T]$ be an optimal control of the OCP \eqref{eq:ocp_cls}, then the control function $\tilde{u}^*(\tau), \tau \in [0,\tilde{\Tctrl}]$ obtained via Eq. \eqref{eq:optimal_control} is optimal
	for \eqref{eq:ocp_time_freezing}.
\end{theorem}
\textbf{PROOF}.
For a fixed control function $u(\tau)$ and $s(\tau)$ such that $t(\tilde{\Tctrl})=\Tctrl$, the time-freezing system \eqref{eq:ocp_time_freezing_dcs_1}-
\eqref{eq:ocp_time_freezing_dcs_comp2} and the CLS \eqref{eq:cls_friction} with $u(t(\tau))$ are equivalent in the sense of Theorem \ref{th:solution_relationship_friction}.
Thus, a feasible $y(\tau)$ in \eqref{eq:ocp_time_freezing} results in a feasible $x(t)$ in \eqref{eq:ocp_cls}. 
Due to equation \eqref{eq:quadrature_tf}, both OCPs have the same objective value.
Consequently, given a $\tilde{u}(\tau) = \tilde{u}^*(\tau) + \delta \tilde{u}(\tau)$ that improves the objective \eqref{eq:ocp_time_freezing_objective}, the corresponding $\tilde{u}(t(\tau))$ would also improve \eqref{eq:ocp_cls_objective}. 
Conversely, for every modified $u(t(\tau)) = u^*(t(\tau)) + \delta u(t(\tau))$ that improves the objective \eqref{eq:ocp_cls_objective}, we can construct an appropriate control function $u(\tau)$ via \eqref{eq:optimal_control} that improves \eqref{eq:ocp_time_freezing_objective}.
Thus, $u^*(t) = \tilde{u}^*(t(\tau))$ is optimal for \eqref{eq:ocp_cls}.
The converse is proved by similar arguments. 
\qed
\subsection{Discrete-time OCP with the time-freezing system}
In principle, one can discretize the OCP \eqref{eq:ocp_cls} by using any time-stepping integration method for CLS \cite{Acary2008,Brogliato2016,Stewart2000} e.g., the Stewart-Trinkle method \cite{Stewart1996}. 
Such an approach for direct optimal control was used in \cite{Posa2014}. 
As discussed in Section \ref{sec:introduction}, standard time-stepping methods for CLS with friction \eqref{eq:cls_friction} have at best first-order accuracy \cite{Acary2008,Stewart2000}.
Moreover, the numerical sensitivities obtained from such a discretization are always wrong and the NLP solvers converge to spurious solutions \cite{Stewart2010,Zhong2022}.
Therefore, for a moderately accurate solution usually a large computational effort is needed.
Standard time-stepping methods for PSS encounter the same difficulties as methods for CLS.
These fundamental limitations motivate the derivation of the OCP formulation with a time-freezing system \eqref{eq:ocp_time_freezing} since we can use the recently introduced FESD method that overcomes these difficulties \cite{Nurkanovic2022}. 
It discretizes dynamic complementarity systems that are equivalent to PSS (such as Eq. \eqref{eq:time_freezing_dcs}) and ensures automatic switch detection, higher-order integration accuracy, and correct numerical sensitivities \cite{Nurkanovic2022b,Nurkanovic2022}. 
In conclusion, the fundamental limitations of standard direct optimal control methods are resolved by combining time-freezing and FESD. 
This enables one to find a more accurate solution approximation for the continuous-time OCP~\eqref{eq:ocp_cls} by solving~\eqref{eq:ocp_time_freezing}.

We proceed by introducing the discrete-time version of the OCP \eqref{eq:ocp_time_freezing} with a multiple shooting-type discretization~\cite{Bock1984}.
The numerical time horizon $[0,\tilde{\Tctrl}]$ is split into $\Nctrl$ control intervals $[\tau_k,\tau_{k+1}]$ of equal length \cite{Nurkanovic2022}.
The controls are assumed to be constant over every interval, i.e., $u(\tau) = u_k, \tau \in [\tau_k,\tau_{k+1}], k = 0,\ldots, \Nctrl-1$, and $y_k = (x_k,t_k)\in \R^{n_y}$ is the discrete-time approximation of the time-freezing state, i.e., $x_k \approx x(\tau_k)$ $t_k \approx t(\tau_k)$. 
The vectors $z_k$ collect all algebraic and internal integration variables for the $k-$th control interval. 
The vector ${w} \coloneqq (y_0,z_0,u_0,s_0,\ldots,y_{\Nctrl\!-\!1},z_{\Nctrl\!-\!1},u_{\Nctrl\!-\!1},s_{\Nctrl\!-\!1},y_{\Nctrl})$ groups all optimization variables.

Our goal is to have an equidistant control grid, as this is typically required in feedback control applications.
It is important to note that, due to intervals with frozen physical time evolution ($t'=0$), an equidistant grid in numerical time $\{\tau_0,\ldots,\tau_{\Nctrl}\}$ does not imply an equidistant grid in physical time $\{t_0,\ldots,t_{\Nctrl}\}$. 
To address this issue, we propose to use a piecewise constant discretization of the time-transformation variable $s(\tau)$, i.e., we have $s_k\in \R, k = 0,\ldots,\Nctrl-1$.
Additionally, we introduce the constraints $t_k = k\frac{\Tctrl}{\Nctrl}$, $k = 0,\ldots,\Nctrl$, cf. Eq. \eqref{eq:ocp_time_freezing_discrete_clock} below.
It is worth noting that for $k = N$, we have the discrete-time versions of the terminal clock constraint \eqref{eq:ocp_time_freezing_temrinal_time}.
The steps above result in an equidistant control discretization grid in physical time, i.e., $u(t) = u_k$ for $t\in [t_k,t_{k+1}]$ with $t_0 = 0$ and $t_k = t_{k-1} + \Tctrl/\Nctrl$. 
 This is further illustrated in Example \ref{ex:illustrate_ocp}.

The discretization of \eqref{eq:ocp_time_freezing} reads as:
\begin{subequations}\label{eq:ocp_time_freezing_discrete}
	\begin{align}
		\underset{{w}}{\min} \;   & \Psi(x_{\Nctrl}) \\
		\textrm{s.t.} \quad 
		&x_0 = \bar{x}_0,\\
		&y_{k +1}\! = \Phi_f(y_k,z_k,u_k,s_k), \, k = 0,\! \ldots,\! \Nctrl\!-\!1,	\label{eq:ocp_time_freezing_discrete_state} \\
		&0 \! =\! \Phi_{\mathrm{int}}\!(y_k,z_k,u_k), \, k = 0,\! \ldots,\! \Nctrl\!-\!1,\label{eq:ocp_time_freezing_discrete_int}\\
		&0\! \leq\! \Phi_{\mathrm{c,1}}\!(z_k)\! \perp \! \Phi_{\mathrm{c,2}}(z_k)\! \geq\! 0,\, k = 0,\! \ldots,\! \Nctrl\!-\!1, \label{eq:ocp_time_freezing_discrete_comp}\\
		&t_k =  k \frac{\Tctrl}{\Nctrl},\; k = 0, \ldots, \Nctrl, \label{eq:ocp_time_freezing_discrete_clock}\\
		&1 \leq s_k \leq \bar{s},\; k = 0, \ldots, \Nctrl-1, \label{eq:ocp_time_freezing_discrete_clock_bounds}\\
		&0 \leq g(x_k,u_k), \; k = 0, \ldots, \Nctrl\!-\!1,\\
		&0 \leq r(x_{\Nctrl}).
	\end{align}
\end{subequations}
It is common in direct optimal control to write discretization method equations in a compact discrete-time system manner \cite[Chapter 8]{Rawlings2017}, as we do here in Eq.~\eqref{eq:ocp_time_freezing_discrete_state}-\eqref{eq:ocp_time_freezing_discrete_comp}.
The function $\Phi_f: \R^{n_y} \times \R^{n_z} \times \R^{n_u} \times \R \to \R^{n_x}$ is the discrete-time state transition map which approximates $y(\tau)$. 
The function $\Phi_{\mathrm{int}}: \R^{n_y} \times \R^{n_z} \times \R^{n_u} \to \R^{n_\Phi}$ collects all internal computations of the underlying integration scheme. 
The constraints \eqref{eq:ocp_time_freezing_discrete_comp} arise from the discretization of the complementarity conditions~\eqref{eq:ocp_time_freezing_dcs_comp1}-\eqref{eq:ocp_time_freezing_dcs_comp2}.
These functions are obtained via the FESD discretization, for more details see \cite[Section 3.2.4]{Nurkanovic2022}. 
The constraint \eqref{eq:ocp_time_freezing_discrete_clock_bounds} bounds $s_k$, where $\bar{s}$ is its maximal value that has to be sufficiently large to ensure feasibility of~\eqref{eq:ocp_time_freezing_discrete_clock}.

The FESD method for PSS is implemented in the CasADi \cite{Andersson2019} based open-source tool NOSNOC \cite{Nurkanovic2022c,Nurkanovic2022b}.
Note that the NLP \eqref{eq:ocp_time_freezing_discrete} is a mathematical program with complementarity constraints. 
\color{black}
They are degenerate nonsmooth NLP which are solved in NOSNOC with a homotopy approach, cf. \cite[Section IV.B]{Nurkanovic2022b}. 
The advantage of the homotopy approach is that only a (finite) sequence of related, but smooth NLP is solved. 
Under some regularity assumptions, the solution of the last NLP is a solution of the initial nonsmooth NLP \cite{Anitescu2007,Hall2022}. 
The main drawback of the homotopy approach is that it sometimes requires some tuning of the homotopy parameters.

\begin{figure}[t]
	\vspace{-0.2cm}
	{\includegraphics[scale=0.70]{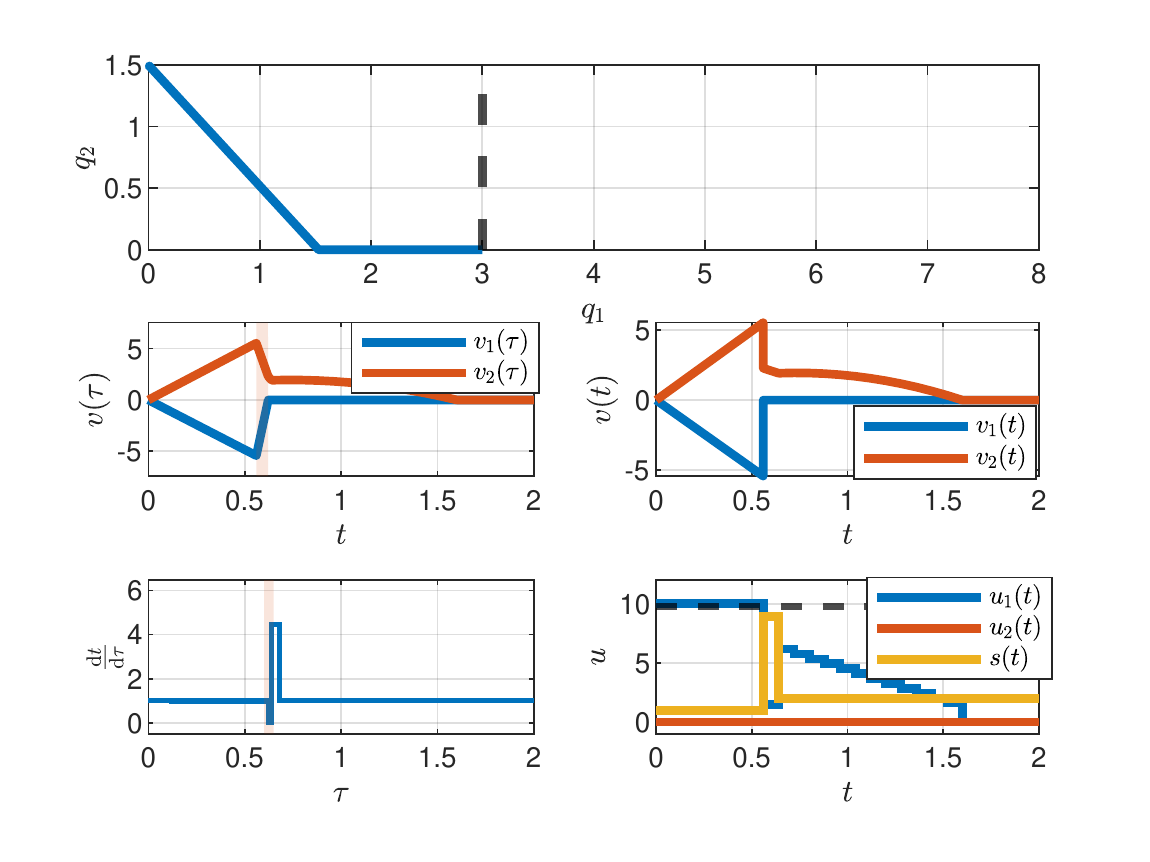}}
	\vspace{-0.7cm}
	\caption{Solution to the guiding optimal control example.}
	\label{fig:ocp_solution}
\end{figure}

\begin{example}(OCP Example)\label{ex:illustrate_ocp}
We solve an OCP of the form of \eqref{eq:ocp_cls} with our guiding example.
The initial value is unchanged, i.e., $y_0 = (0,1,0,0,0)$. 
The particle should reach at $t(\tilde{\Tctrl}) = 2$, with $\tilde{\Tctrl}=2$, the position $q(T) = (3,0)$ with zero terminal velocity $v(T) = (0,0)$. 
We bound the horizontal thrust force $|u_1| \leq 10$ and set for simplicity $u_2 = 0$. 
The ball should reach the goal with minimum control effort, which is modeled with the stage cost $L(x,u) = u_1^2$.
We take $\Nctrl = 20$ control intervals and discretize the equivalent time-freezing OCP with a third-order FESD-Radau II-A scheme with three integration steps on every control interval~\cite{Nurkanovic2022b}. 
The solution is depicted in Fig.~\ref{fig:ocp_solution}.
We can see that maximum force is applied before the impact since there is still no friction and the motion is {cheaper}. 
After the impact, a smaller control force is applied just to reach the target. 
Note that the $s_k$ (yellow line in the bottom right plot) is higher during the control interval when the state jump happens (to catch up the frozen time) and $s_k=2$ during contact to compensate for the slow down due to $\gamma(x,u)$. 
The resulting speed of time is always one (bottom left plot), except when the state jump happens where a speed-up is needed to compensate for the frozen time. 
This ensures an equidistant control grid (as intended with the constraint \eqref{eq:ocp_time_freezing_discrete_clock}) and $t(\tilde{\Tctrl}) = 2$ as desired.
\end{example}

\section{Numerical optimal control of a jumping robot}\label{sec:numerical_examples}
\begin{figure*}[t!]
	\centering
	\centering
	\vspace{0.2cm}
	{\includegraphics[scale=0.54]{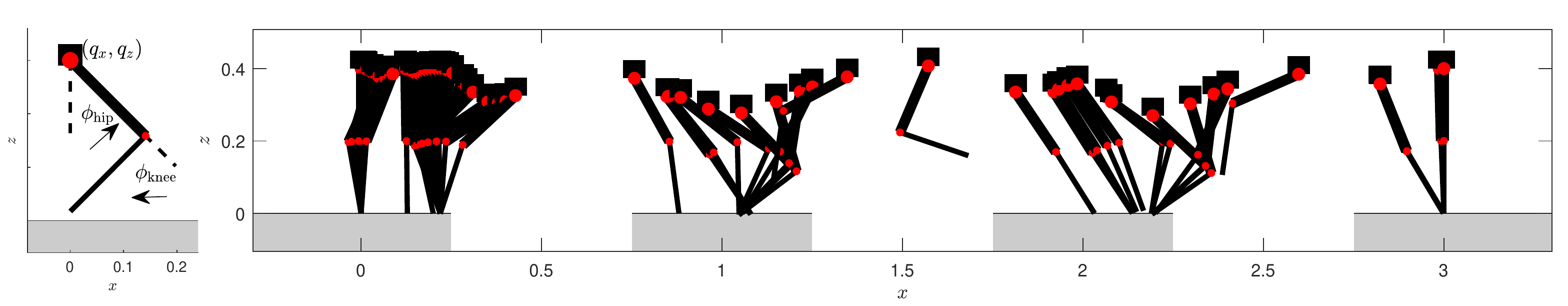}}
	\vspace{-0.4cm}
	\caption{Illustration of the robot kinematics (left), several frames of the solution of the discretized OCP (right).}
	\label{fig:hopper_ocp}
	\centering
\end{figure*}

We consider a hopping robot that must jump over three holes to reach a desired target. 
This example showcases all theoretical developments of Sections \ref{sec:time_freezing}, \ref{sec:frictional_impact}, and the use of numerical methods described in Section~\ref{sec:computational_considerations}.
Thereby, an OCP formulation for synthesizing dynamic motions of the single-legged 2D robot \textit{Capler} \cite{Carius2018,Hwangbo2018} is derived.
The robot is described by four degrees of freedom $q = (q_x,q_z,\phi_\mathrm{knee},\phi_\mathrm{hip})$.
Here, $(q_x,q_z)$ are the coordinates of the robot's base at the hip and $\phi_\mathrm{knee}$, $\phi_\mathrm{hip}$ are the angles of the hip and knee, respectively, cf. left plot in Fig. \ref{fig:hopper_ocp}.
It is actuated by two direct-drive motors at the hip and knee joints. 
The robot's dynamics are compactly described by the CLS in the form of \eqref{eq:cls_friction}.
The torques of the two motors $u(t)=(u_{\mathrm{knee}}(t),u_{\mathrm{hip}}(t))$ are the control variables.
A detailed derivation of the model equations and all parameters for the robot can be found in~\cite[Appendix A]{Gehring2011}. 

Denote by $p_{\mathrm{foot}}(q) = (p_{\mathrm{foot},x}(q),p_{\mathrm{foot},z}(q))$ and  $p_{\mathrm{knee}}(q) = (p_{\mathrm{knee},x}(q),p_{\mathrm{knee},z}(q))$ the kinematic position of the robot's foot and knee, respectively. For the unilateral constraint function we take $f_c(q) = p_{\mathrm{foot},z}(q)$. 
For a planar robot, we need just one tangent, i.e., $b(q) = \nabla_q p_{\mathrm{foot},x}(q)$ and the friction model is exact.
The coefficient of friction is $\mu = 0.8$ and the auxiliary ODE constant is $\anormal = 200$. 

The objective of the OCP is to minimize the integral of the squared control torques, i.e., we have the stage cost $L(x,u) = u(\tau)^\top u(\tau)$.
The robots should reach a given target position $q_{\mathrm{target}}=(3,0.4,0,0)$ starting from the initial position 
$q_{0}=(0,0.4,0,0)$ with zero velocity $v_{0}=\mathbf{0}_{4,1}$.
The initial value is $y_0 = (q_0,v_0,0)$.

The prediction horizon is  $\Tctrl = 2.5$ s.
We add the following constraints on the states and kinematic positions:
\begin{align*}
	-0.05 e&\leq (q_x(t),p_{\mathrm{foot},x}(q(t)),p_{\mathrm{knee},x}(q(t))),\\
	 	0.2&\leq q_z(t) \leq 0.55,\\
	 -\frac{3\pi}{8}&\leq\phi_\mathrm{hip}(t)\leq \frac{3\pi}{8},\\
	 -\frac{\pi}{2} &\leq \phi_\mathrm{knee}(t) \leq \frac{\pi}{2},\\
	 	 0.05&\leq p_{\mathrm{knee},z}(q(t)),\\
	 -0.005&\leq 	p_{\mathrm{foot},z}(q(t)) \leq 0.2,\; t\in [0,T].
\end{align*}
Their goal is twofold. On one hand, they should avoid unnatural and too extensive bending of the joints. On the other hand, they serve as \textit{guiding constraints} during the early phases of the homotopy procedure. In the early iterations, the physics is relaxed and we want to prohibit the optimizer to go to undesired regions.
The control bounds read as 
\begin{align*}
	-60 e &\leq u(t) \leq 60 e,\; t\in [0,T].
\end{align*}

On the way to the target, the robot must overcome three holes in the ground. Instead of using very complicated expressions for $f_c(q)$ we model the holes as regions that the robot should not enter. This is achieved by constraints inside the OCP requiring that $p_{\mathrm{foot}}$ is outside $n_e = 3$ ellipsoids:
\begin{align*}
\Big(\frac{p_{\mathrm{foot},x}-x_{\mathrm{c},k}}{a_k}\Big)^2 + \Big(\frac{p_{\mathrm{foot},z}-z_{\mathrm{c},k}}{b_k}\Big)^2 \geq 1, k = 1,\dots, n_e.
\end{align*}
By appropriately picking $a_k,b_k,x_{\mathrm{c},k}$ and $z_{\mathrm{c},k}$ the desired shapes are trivially selected. In our example, we pick $z_{\mathrm{c},k} =0$, $a_k = 0.5$ (width of the hole), $b_k = 0.1$ (kept low, should not enforce unnecessarily high jumps). For the centers of the holes, we pick $x_{\mathrm{c},1} = 0.5$, $x_{\mathrm{c},2} = 1.5$, $x_{\mathrm{c},3} = 2.5$. We collect all path constraints (for the holes, on the kinematics, control, and state bounds) into the function $g(x,u) \geq 0$. 
\begin{remark}
Note that the constraints $g(x,u) \geq 0$ cannot become active if the corresponding normal velocity is nonzero, as opposed to activating a constraint $f_c(q)\geq 0$, since no state jump law is associated with path constraints in an OCP. This is one of the main differences between constraints that are part of the dynamics (equipped with a state jump law) and path constraints in the OCP.
\end{remark}

\begin{figure}[t]
	\centering
	\vspace{0.1cm}
	\centering
	{\includegraphics[scale=0.65]{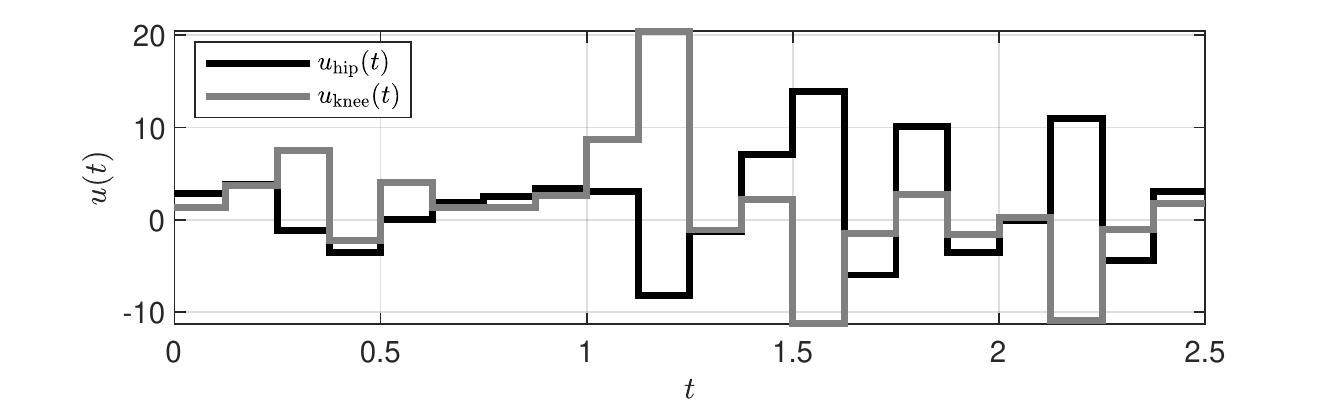}}
	\vspace{-0.6cm}
	\caption{The optimal control input $u(t)$ in physical time $t$ obtained by solving the discretization of the optimal control problem.}
	\label{fig:ocp_controls}
\end{figure}

We have now all ingredients to formulate an OCP of the form of Eq. \eqref{eq:ocp_cls}. NOSNOC automatically reformulates the CLS into a time-freezing system and discretizes the OCP, such that we obtain a discrete-time problem of the form of \eqref{eq:ocp_time_freezing}.
The resulting mathematical program with complementarity constraints is solved in a homotopy procedure with IPOPT \cite{Waechter2006} equipped with the MA57 linear solver \cite{HSL}.
The source code for this example is available in NOSNOC's repository \cite{Nurkanovic2022c}.
The OCP is discretized with a FESD Radau-IIA scheme of order 5 \cite{Nurkanovic2022}. 
We consider $\Nctrl=20$  control intervals with $3$ intermediate integration steps on every interval.

For the initialization of the differential states, we take $y_0$ at every discretization node. All discrete-time control variables are initialized with zero. Hence, no information about the order, number, or timing of the nonsmooth transitions and jumps is provided. Treating the contact dynamics directly in the OCP and thus implicitly discovering all nonsmooth transitions is in the robotics community called \textit{contact implicit optimization} \cite{Carius2018,Tassa2012}. The results of the optimization are shown in the right plot of Fig.~\ref{fig:hopper_ocp}. The approach finds an intuitive dynamic movement by solving only smooth NLP, without providing any hints about the order and number of nonsmooth transitions. The optimal torques are depicted in Fig.~\ref{fig:ocp_controls}.

\section{Conclusion and outlook}\label{sec:conclusions}
This paper introduced a novel time-freezing reformulation for transforming complementarity Lagrangian systems (CLS) with inelastic impacts and friction into piecewise smooth systems. 
This reduces the level of nonsmoothness significantly. 
We prove solution equivalence under mild conditions and derive constructive ways to select the auxiliary ODE.
Moreover, we show that we can obtain a solution to an optimal control problem with a CLS by solving the simpler OCP with time-freezing systems.
We derive a reformulation of the time-freezing Filippov system into a dynamic complementarity system that allows the use of the high-accuracy Finite Elements with Switch Detection (FESD) method for numerical optimal control \cite{Nurkanovic2022}.
To the best of our knowledge, this is the first reformulation that enables one to treat CLS with inelastic impacts as Filippov systems. 
The practicality of the discussed methods is demonstrated in an OCP considering a one-legged robot with frictional impact.  
All methods from this paper are implemented in the open-source package NOSNOC \cite{Nurkanovic2022b}.

In further work, we aim to extend the ideas to multiple and simultaneous impacts, e.g., based on Moreau's impact law in second-order sweeping processes \cite{Brogliato2016}. 
Moreover, to make this approach more practical, good initialization strategies and more sophisticated homotopy procedures for the mathematical programs with complementarity constraints would be useful.

\begin{ack}                               
We thank Jan Carius from ETH Z\"urich in Switzerland, for providing details for the \textit{Capler} robot.
\end{ack}

\appendix

\section{Proof of Theorem \ref{th:solution_relationship}}\label{app:proof_of_th1}
\textbf{PROOF.}
The idea of the proof is to consider the different modes, in which the CLS and time-freezing system can be and to compare the solutions to establish the result of the theorem.
A solution of the initial value problem given by the time-freezing system in Eq. \eqref{eq:time_freezing_system} with $y(0) = y_0$ is denoted by $y_{\mathrm{sol}}(\tau;y_0)$ for $\tau \in [0,\hat{\tau}]$. Similarly, for the CLS in Eq. \eqref{eq:cls} and $x(0) = x_0 $ for $t(\tau) \in [0,t(\hat{\tau})]$ we use $x_{\mathrm{sol}}(t(\tau);x_0)$. We must distinguish all possible cases, hence we split the proof into several parts. 

\vspace{-0.1cm}
\textit{Part I (Unconstrained case).}
Regard the case $f_c(q(\tau)) >0, \tau \in [0,\hat{\tau}]$. This means that $y\in R_1$, $y' = f_1(y,u) =  ({f}_{\mathrm{ODE}}(x,u),1), \tau \in [0,\hat{\tau}]$. It holds that $t(\tau) = \int_0^{\tau} \dd s = \tau$ and by setting $\hat{\tau} = \tauf$ we have $t(\tauf) = \tf$. 
Note that $Ry' = R {f}_{1}(y,u)$, is equivalent to $ x' = {f}_{\mathrm{ODE}}({x,u})$. Since $[0,\tauf] = [0,\tf]$, this ODE has the same solution as $\dot{x} = f_{\mathrm{ODE}}(x,u)$, therefore relation \eqref{eq:SolutionRelation1_diff} holds for $t \in [0,\tf]$. This means that $f_c({q}(t))>0$ and $\zn(t) = 0$  for  $t \in [0,\tf]$. For the time-freezing system this means that $y(\tau)\notin \Sigma$ for $\tau \in [0,\tauf]$, hence equation \eqref{eq:SolutionRelation1_alg} is also satisfied. 

\textit{Part II (Sliding mode/persistent contact).} 
Regard the case $f_c(q(0)) = 0$ and $n({q}(0))^\top  {v}(0) = 0$, i.e., $y(0) \in  {\Sigma}$. Assume that $\varphi(x(\tau),u(\tau)) <0,\; \tau \in [0,\tauf]$.  This means that $y(\tau)\in \Sigma,\; \tau \in [0,\tauf]$, cf. Section \ref{sec:presistent_contact}. From Theorem \ref{th:presistent_contact} we have that 
$y' =  \gamma(x,u)({f}_{\mathrm{DAE}}(x,u),1)$ for $\tau \in [0,\tauf]$ and $t'(\tau) = \gamma(x,u)>0, \tau  \in [0,\tauf]$, thus  $t(\tauf) = \tf>0$.
On one hand, from $\frac{\dd y(t(\tau))}{\dd \tau} = \frac{\dd y}{\dd t}  \frac{\dd t}{\dd \tau}$ we have that
$R\frac{\dd y}{\dd t} = R\frac{1}{\gamma(x,u)}\gamma(x,u)({f}_{\mathrm{DAE}}(x,u),1) = {f}_{\mathrm{DAE}}(x,u)$. 
On the other hand, for the CLS we have for $f_c({q}(0)) = 0$, $n({q}(0))^\top v(0) = 0$, $z(t)\geq 0,\ t \in [0,\tf]$. Consequently, the CLS reduces to the ODE $\dot{x} = f_{\mathrm{DAE}}(x,u)$. Similar to part I, we conclude that \eqref{eq:SolutionRelation1_diff} holds. Since $f_c(q(t(\tau))) =0$ and $n(q((\tau)))^\top v((\tau)) =0$ for $\tau \in [0,\tauf]$, the relation for $\zn(t(\tau))$ in Eq. \eqref{eq:SolutionRelation1_alg} follows from Eq. \eqref{eq:contact_LCP_solution_map} and \eqref{eq:cls_index_reduced_dae}.

\textit{Part III (Leaving sliding mode).} 
Now we consider a similar scenario as in part II, with  $y_0 \in  {\Sigma},\; \varphi(x(\tau),u(\tau))\leq 0$ with $\tau \in \left[0,\taue\right),\taue < \tauf$ (sliding mode) and $\varphi(x(\tau),u(\tau))>0$ for $\tau \in [\taue,\tauf]$ (leaving sliding mode).
Relations \eqref{eq:SolutionRelation1_diff} and \eqref{eq:SolutionRelation1_alg} hold for $\tau \in \left[0,\taue\right)$ by the same arguments as in part II.
For $\tau \geq \taue$, following the arguments in Section \ref{sec:contact_breaking}, $y(\tau)$ leaves $\Sigma$ and $y(\tau) \in R_1$ for $\tau \in \left[0,\taue\right)$. We can apply the arguments of part I and establish the result of the theorem.

\textit{Part IV (State jump).} This part regards the case of $\taus \in [0,\tauf]$, i.e., $\ts \in [0,\tf]$.
For $\tau \in \left[0,\taus\right)$ and $t\in \left[0,t(\taus^-)\right)$ we can apply Part I of the proof by simply setting $\hat{\tau} = \taus$ and deduce that \eqref{eq:SolutionRelation1_diff} and \eqref{eq:SolutionRelation1_alg} hold.
For $\tau = \taus$ we have $f_c(q(\taus)) = 0$ and $n({q}(\taus))^\top  {v}(\taus)<0$. Consequently, $y\in R_2$ and $y' =  {f}_{\mathrm{aux,n}}(y)$. The assumption 
$f_c(q(\tau_{\mathrm{f}})) \geq 0$ and $n({q}(t_{\mathrm{f}}))^\top {v}(t_{\mathrm{f}})\geq0$ ensures that the time evolution of $y'(\tau) = f_{\mathrm{aux,n}}(y(\tau))$ is finished in $[\taus,\tauf]$, i.e., $\taur \leq \tauf$.
From the proof of Proposition \ref{prop:aux_forming_ode} we know that by construction $ {q}(\tau) =  {q}(\taus) \eqqcolon  \qs,\ \tau \in [\taus,\taur]$. 
Consequently, $f_c(q(\tau)) = 0 ,\ \tau \in [\taus,\taur]$. 
For $v(\tau)$, from \eqref{eq:aux_dyn_example} we obtain that: 
\begin{align}\label{eq:v_solution_virtual_time}
 {v}(\taur)&=  {v}(\taus) + \int_{\taus}^{\taur} M(q(\tau))^{-1} n(q(\tau))\anormal \dd \tau.
\end{align}
Multiplying both sides with $ n( \qs)^\top$ from the left and noting that $M^{-1}(q(\tau))n({q}(\tau))$ is constant since $q(\tau) = \qs,\ \tau \in [\taus,\taur]$, we have
\begin{align*}
&\underbrace{n(\qs)^\top \! v(\taur)}_{=0}\!-\! n( \qs)^\top\! v(\taus)\! =\!  \underbrace{n( \qs)^\top\! M( \qs)^{-1}\!n( \qs)}_{=D(\qs)} \underbrace{\int_{\taus}^{\taur} \!\! \anormal \dd \tau}_{\eqqcolon \Lambda_1},\\
&\Lambda_1= -\frac{ n(  \qs)^\top v(\taus)}{D( \qs)}>0.
\end{align*}
Next, we look at the post-impact states of the CLS and compare it to the solution of the time-freezing system.
Since in CLS, $v(t)$ is a function of bounded variation \cite{Brogliato2016}, we have that $q(t)$ is a continuous function. Thus, $q(\ts^+) =	q(\ts^-) = q(\ts)$. Furthermore, notice that $ q(\ts)= \qs$ which implies $n(q(\ts)) ^\top M^{-1}(q(\ts)) n(q(\ts)) = D( \qs)$.
Examining,
\begin{align}\label{eq:v_solution_physical_time}
\begin{split}
v(\ts^+) &=	v(\ts^-) +  \underbrace{\lim_{\substack{\epsilon \to 0\\ \epsilon>0}}\int_{\ts-\epsilon}^{\ts+\epsilon}f_\mathrm{v}(q(t),v(t)) \dd t }_{=0}\\
&+ \lim_{\substack{\epsilon \to 0\\ \epsilon>0}}\int_{\ts-\epsilon}^{\ts+\epsilon}{M(q(t))^{-1} n(q(t))}\zn(t)\dd t,
\end{split}
\end{align}
and multiplying both sides with $n(q_\mathrm{s})^\top$ from the left, introducing $\Lambda_2 \coloneqq  \lim_{\substack{\epsilon \to 0\\ \epsilon>0}}\int_{\ts-\epsilon}^{\ts+\epsilon}\zn(t)\dd t$,   we conclude that 
\begin{align}\label{eq:v_intergral_physical_time}
	\Lambda_2 = -\frac{n(q(\ts))^\top v(\ts^-)}{D( \qs)} = \Lambda_1.
\end{align}
By comparing \eqref{eq:v_solution_virtual_time} and \eqref{eq:v_solution_physical_time}, due to the last relation we conclude that $ {v}(\taur) = v(\ts^+) \eqqcolon   \vs$. 
Furthermore, by Proposition \ref{prop:aux_forming_ode} we have $f_c(q(\taur)) = 0$ and $n({q}(\taur))^\top  {v}(\taur) =0$. 
Since $t'(\tau) = 0$ with $\tau \in [\taus,\taur]$, it follows that $t(\taur) = t(\taus) = t_s$.
Consequently,
\begin{subequations}\label{eq:consistent_init_proof}
	\begin{align}
		f_c(q(\ts)) &= f_c(q(\taur)) =0,\\
		n(q(\ts))^\top {v}(\ts^+) &=  n(q(\taur))^\top v(\taur)=0.
	\end{align}
\end{subequations}
Let $y_{\mathrm{s}} \coloneqq ( \qs, \vs,\ts)$. Note that $y_{\mathrm{sol}}(\tau - \taur, y_{\mathrm{s}}) = y(\tau,y_0)$ for $\tau \in [\taur, \tauf]$. Likewise, $x_{\mathrm{sol}}(t - \ts, x_{\mathrm{s}}) = x(t,x_0)$ for $t \in \left(\ts, \tf\right]$, with $x_{\mathrm{s}} = R y_{\mathrm{s}}$.  The two initial value problems are initialized with the same initial condition. Since \eqref{eq:consistent_init_proof} holds we can apply Theorem \ref{th:presistent_contact} for $y\in \Sigma$. Therefore, by using the arguments of parts II or III (depending on $\varphi(x(\tau),u(\tau))$), we deduce that \eqref{eq:SolutionRelation1_diff} and \eqref{eq:SolutionRelation1_alg} hold on $[\taur, \tauf]$.
Additionally, for $\tau \in (\taus,\taur)$ we have $t=\ts$ and Eq. \eqref{eq:SolutionRelation2} follows directly from \eqref{eq:v_intergral_physical_time}.

\textit{Part V (Summary).}
Parts I-IV cover all possible modes of the CLS and the time-freezing system: evolution according to $f_{\mathrm{ODE}}$ (Part I), evolution on $\Sigma$ according to $f_{\mathrm{DAE}}$ without leaving it (Part II), leaving $\Sigma$ and continuing to evolve according to $f_{\mathrm{ODE}}$ (Part III), and the state jump (Part IV). To regard any other possible sequence of mode on $[0,\tauf]$, the time interval is simply split into sub-intervals with the different mode, and we apply subsequently the arguments from Parts I-IV to verify that \eqref{eq:SolutionRelation1_diff} and \eqref{eq:SolutionRelation1_alg} hold for $t\neq \ts$ and \eqref{eq:SolutionRelation2} for $t=\ts$. This completes the proof.\qed

\section{Proof of Theorem \ref{th:slip_stick_theorem}}\label{app:proof_of_th2}
\textbf{PROOF.}
For the first part of the proof, we assume that $b(q)^\top v > 0$. This means that $y \notin R_3$ and it follows that $\theta_3 = 0$. 
From the conditions $c_2(y) = 0$ and $\theta_1 + \theta_2 = 1$ we can compute $\theta_1$ and $\theta_2$. By using the fact that $\nabla n(q)^\top M(q)^{-1} b(q) = 0$ and following similar lines as in the proof of Theorem \ref{th:presistent_contact} we can compute that 
$y' =(v,f_{\mathrm{v}}(q,v,u)-M(q)^{-1} D(q)^{-1}\varphi(x,u)  (n(q)-\mu b(q))$. By similar reasoning for $b(q)^\top v < 0$ we obtain that $y' =(v,f_{\mathrm{v}}(q,v,u)-M(q)^{-1} D(q)^{-1}\varphi(x,u)(n(q)+\mu b(q))$. Therefore, it holds for $\vt \neq 0$ that ${y'=\gamma(x,u)(\fslip(x,u),1)}$, whereby the $\mathrm{sign}(\cdot)$ in the r.h.s. of $\fslip(x,u)$ accounts for the sign of~$\vt$.

In the second part, we have $ \vt = b(q)^\top v =0$. Together with the assumption that $y\in \Sigma$ it follows that $y\in \partial R_i$, for all $i =1,2,3$, hence no $\theta_i$ can be set to be zero a priori.
Since the vectors $b(q)$ and $n(q)$ are orthogonal in the kinetic metric, it can be seen that the matrix $\tilde{D}(q)$
is a diagonal matrix. We denote the first and second entries on its diagonal by ${D}(q)$ (cf. \eqref{eq:contact_LCP_functions_D}) and ${D}_{\mathrm{t}}(q)$, respectively. Similarly, the first and second components of the vector $\tilde{\varphi}(x,u)$ are denoted by $\varphi(x,u)$ and $\varphi_{\mathrm{t}}(x,u)$, respectively. This, we obtain that
 $\zn = -D(q)^{-1}\varphi(x,u)$ and $\zt = D_{\mathrm{t}}(q)^{-1} \varphi_{\mathrm{t}}(x,u)$ holds.

From the condition $\frac{\dd c_2(y)}{\dd \tau} = 0$ and the definition of the time-freezing system we compute that
\begin{align*}
	\varphi(x,u) \theta_1 + D(q)\anormal(\theta_2+\theta_3) = 0.
\end{align*}
Using $\sum_{i=1}^{3} \theta_i = 1$ and the last equation we can compute that:
\begin{subequations}\label{eq:stick_proof_normal}
\begin{align}
	\theta_1 &= \frac{D(q)\anormal}{D(q)\anormal-\varphi(x,u)} = \gamma(x,u), \label{eq:stick_proof_normal_theta}\\
	\theta_2+\theta_3&= \frac{-\varphi(x,u)}{D(q)\anormal-\varphi(x,u)} \label{eq:stick_proof_normal_sum_theta}
\end{align}
\end{subequations}
Next, we use the condition $\frac{\dd c_3(y)}{\dd \tau}= \frac{\dd (b(q)^\top v)}{\dd \tau} = 0$ and compute:
\begin{align*}
	\varphi_{\mathrm{t}}(x,u) \theta_1  + \mu {D}_{\mathrm{t}}(q)\anormal (-\theta_2 +\theta_3 ) = 0,
\end{align*}
and by using \eqref{eq:stick_proof_normal_theta} we establish the relation:
\begin{align}\label{eq:stick_proof_tan1_sum_theta}
-\theta_2 +\theta_3&=  -\gamma(x,u) \frac{ \varphi_{\mathrm{t}}(x,u)}{\mu D_{\mathrm{t}}(q)\anormal}.
\end{align}
Next, we compute the sliding mode vector field of the time-freezing system $y' = \sum_{i=1}^{3}\theta_i f_i $. By rearranging the terms we obtain
\begin{align*}
\begin{split}
y' &= \theta_1\begin{bmatrix}
v \\ f_v \\ 1 
\end{bmatrix}
+(\theta_2+\theta_3)\begin{bmatrix}
	\mathbf{0}_{n_q,1} \\ M(q)^{-1}n(q) \anormal\\ 0 
\end{bmatrix} \\
&+(-\theta_2+\theta_3)\begin{bmatrix}
	\mathbf{0}_{n_q,1} \\ M(q)^{-1}b(q) \mu  \anormal\\ 0 
\end{bmatrix} 
\end{split}
\end{align*}
Now, by multiplying the second term by $D(q)^{-1} D(q)$, plugging in the expressions for the sums  \eqref{eq:stick_proof_normal_sum_theta} and \eqref{eq:stick_proof_tan1_sum_theta} and by comparing it to the r.h.s. of  Eq. \eqref{eq:cls_stick} we can conclude that $y' = \gamma(x,u) (\fstick(x,u),1)$. This completes the proof. \qed

\bibliographystyle{plain}

\end{document}